\documentclass[a4paper,twoside,11pt]{article} % to use for all kind of papers but AACA

\usepackage[english]{babel} %
\usepackage{lmodern}
\usepackage[T1]{fontenc}
\usepackage[latin9]{inputenc}	% Così le "è" vengono trasformate in \`e etc. etc.
\newlength{\defaultparindent}
\usepackage{latexsym}
\usepackage{amsmath}
\usepackage{amsfonts}
\usepackage{amsthm}
\usepackage{amssymb}
\usepackage{bbold}
\usepackage{multirow}
\usepackage{hyperref}
\hypersetup{colorlinks=true,citecolor=green,linkcolor= blue}
\usepackage{todonotes} % more appealing than \marginpar ??
\usepackage{soul} % to get overstriked text

\usepackage{algorithm} % to write algorithms, suggested by Renato on May 11, 2017
\usepackage{algpseudocode}

\usepackage[arXiv]{optional} % Possible values of the option fixing the format: x, std, arXiv, JMP, JOPA, AACA, TCS and also note control: margin_notes, final_notes
\def\mynote{\todo} % \marginpar or \todo
\opt{AACA}{% only for submission to AACA -> \cal font is not defined....
\def\cal{\mathcal}
}

\opt{x,AACA}{% in all cases - but arXiv, std, JMP & JOPA - include my file of standard defs: mbh_defs_TeX
\input{mbh_defs_TeX}% My standard TeX definitions, original in Dictionaries & Gloassaries folder
}

\opt{arXiv,std,JMP,JOPA,TCS}{% only for arXiv, std, JMP, JOPA & TCS we need the definitions included in mbh_defs_TeX (that contains original versions)

\newtheorem{MS_theorem}{Theorem}
\newtheorem{MS_lemma}{Lemma}
\newtheorem{MS_Proposition}{Proposition}
\newtheorem{MS_Corollary}[MS_Proposition]{Corollary}

\def\myconjugate#1{\overline{#1}} % this way conjugate definitions can be easily changed (\bar did not work fine on multiple chars)

\def\eg{e.g.\ }
\newcommand{\R}{\ensuremath{\mathbb{R}}} % good; MB XII 2005; needs \usepackage{bbold}
\newcommand{\Identity}{\ensuremath{\mathbb{1}}} % good; MB XII 2005; needs \usepackage{bbold}
\def\my_span#1{\mbox{Span}\left(#1\right)} % changed from 'span' since it interfered with \multicolumn{} MB X 2009
\def\dotinformula{\;\; \mathrm{.}} % defines space + a full stop (in \rm font) to be placed at the end of a formula
\def\OO#1{\ensuremath{\mbox{O}\!\left(#1\right)}}
\def\SO#1{\ensuremath{\mbox{SO}\!\left(#1\right)}}

\def\O1#1{\ensuremath{\mbox{O}^{#1}(1)}}
\newcommand{\comm}[2]{\ensuremath{\left[ #1, #2 \right]}}
\newcommand{\anticomm}[2]{\ensuremath{\left\{ #1, #2 \right\}}} % or \left[...\right]^+
\newcommand{\myClg}[3]{\ensuremath{{{\cal C}\ell} {\left( #3 \right)}}}	% this is a possible def of myCl that uses only the scalar product g

}

\opt{JMP}{% only for submission to JMP -> double line spacing
}

\def\h_eigen{\eta}
\def\g_eigen{\theta}
\def\mygen{e} % this way the definition of the generators of the algebra can be easily changed
\def\myprimidemp{\mathbb{p}} % this way the definition of the primitive idempotents can be easily changed (beware: not all fonts in math mode have lower characters, e.g. \cal)

\def\SAT{\ensuremath{\mbox{SAT}}}

\def\literal{\ensuremath{\rho}}
\def\mySpinorS{{\mathbb{S}}} % this way the definition of the MLI of spinor spaces S can be easily changed
\def\myFockB{{\cal F}} % this way the definition of the Fock basis can be easily changed
\def\mysetM{\ensuremath{{{\cal M}_n}}} % this way the definition of the set of $2^n$ MTNP of SAT or of Witt basis can be easily changed
\def\mysnqG{\ensuremath{{{\cal N}_n}}} % this way the definition of the semi-neutral quadric Grassmannian can be easily changed
\def\myBooleanS{{\cal{S}}} % this way the definition of the a Boolean SAT problem S can be easily changed
\def\myBooleanT{\mathrm{T}} % this way the definition of the a Boolean T may be easily changed
\def\myBooleanF{\mathrm{F}} % this way the definition of the a Boolean F may be easily changed
\def\myBooleanv{Boolean variable} % this way we switch easily between possible names: logical variable, Boolean variable, literal etc.

\begin{document}

\opt{x,std,arXiv,JMP,JOPA,TCS}{% in all cases - but AACA
\title{{\bf The Boolean SATisfiability Problem and the orthogonal group $\OO{n}$} %\\(temporary title)
	}

\author{\\
	\bf{Marco Budinich}%
%
%\footnote{on leave of absence from: University of Trieste, Trieste, Italy}%
%
\\
%	ICTP and INFN, Trieste, Italy\\
	University of Trieste and INFN, Trieste, Italy\\ % - %\\
	\texttt{mbh@ts.infn.it}\\
%	\texttt{http://www.ts.infn.it/\~{ }mbh/MBHgeneral.html}\\
%
%\\	Very preliminary - restricted circulation (FYEO)
%
%
%	Submitted on March 17, 2016
%	Submitted to: {\em Journal of Mathematical Physics} on September 29, 2017
%		keywords of this submission: Spinor, Satisfiability problem, Clifford algebra
%	Submitted to: {\em Communications in Mathematical Physics} on March 14, 2016.\\
%	Submitted to: {\em Journal of Physics A: Mathematical and Theoretical} \\on October 20, 2017
%	Submitted to: {\em SIAM Journal on Discrete Mathematics} on March 30, 2018
%\\	Re-resubmitted to: {\em Advances in Applied Clifford Algebras} on June 12, 2017
%	Published in: {\em Advances in Applied Clifford Algebras}, 2015\\
%	{\small DOI:10.1007/s00006-015-0547-8}
%	Resubmitted to: {\em Journal of Mathematical Physics} on June 27, 2016\\
%	{\tiny First submission: March 30, 2016}
%	To appear in: {\em Journal of Mathematical Physics}\\
%	{\tiny Submitted on March 30 and June 27, 2016; accepted July 12, 2016}
%	{\em Journal of Mathematical Physics} {\bf 57} (2016) DOI: 10.1063/1.4959531\\
%	\\ May 29, 2017, submitted
%	Submitted to: {\em Theoretical Computer Science} on May 22, 2019
%	Resubmitted to: {\em Theoretical Computer Science} on August 27, 2020
%	Re-resubmitted to: {\em Theoretical Computer Science} on June 8, 2021
	}
\date{ \today }
%\date{April 21, 2017}
%\date{ } % to hide the date this line must be present (believe it or not...)
\maketitle
}

\vspace*{-8mm}% only for this paper to fit all in page 1: this is to move up the Abstract to leave space for keywords, -5/9 is the minimum needed in this case
\begin{abstract}
We explore the relations between the Boolean Satisfiability Problem with $n$ \myBooleanv{}s and the orthogonal group $\OO{n}$. We show that all $2^n$ possible solutions induce involutions of $\R^n$ that lie in the compact, disconnected real manifold of dimension $n (n-1)/2$ of $\OO{n}$. This result in turn gives a new unsatisfiability test within group $\OO{n}$.
\end{abstract}

\opt{x,std,arXiv,JMP,JOPA,TCS}{% in all cases - but AACA
%\vspace*{-1mm}% to move up the keywords - only for this paper to fit all in page 1
{\bf Keywords}: Satisfiability; Clifford algebra; orthogonal group
%\noindent{\bf Keywords:} {Spinors, discrete transformations, Clifford algebra.}
}

\opt{AACA}{% only for AACA
\keywords{Clifford algebra, spinors, Fock basis.}
\maketitle
}

\section{Introduction}
\label{Introduction}
In this paper we explore the relations between the Boolean Satisfiability Problem (\SAT{}) with $n$ \myBooleanv{}s and the orthogonal group $\OO{n}$ of $\R^n$ exploiting the algebraic \SAT{} formulation in Clifford algebra \cite{Budinich_2017} where we proved that a \SAT{} problem is unsatisfiable if and only if it has the maximally symmetric form of the scalars of the algebra.

In Section~\ref{neutral_R_nn} we review some properties of the neutral space $\R^{n,n}$ and of the quadric Grassmannian \mysnqG{} \cite{Porteous_1995}, namely the set of its null subspaces of dimension $n$, within the unifying frame of the Clifford algebra $\myClg{}{}{\R^{n,n}}$. In Section~\ref{SAT_in_Cl(g)} we summarize the relevant parts of \SAT{} formulation in Clifford algebra \cite{Budinich_2017}. In the following three Sections we elaborate on these results formulating \SAT{} problems in a purely geometrical setting. The main result is Theorem~\ref{SAT_in_O(n)} that shows that a given \SAT{} problem is unsatisfiable if and only if the subsets of isometries induced by its clauses cover the orthogonal group $\OO{n}$. To test a \SAT{} problem for unsatisfiability it is thus sufficient to verify if the set of its clauses induce a cover for the compact, disconnected real manifold of dimension $n (n-1)/2$ of $\OO{n}$. This changes the space of solutions to explore from a discrete set of $2^n$ elements to the real manifold of $\OO{n}$. We just remind that, from the computational viewpoint, an unsatisfiability test easily gives an algorithm that actually finds solutions.

For the convenience of the reader we tried to make this paper as elementary and self-contained as possible.

%\newpage
\section{$\R^{n,n}$, its Clifford algebra and $\OO{n}$}
\label{neutral_R_nn}
%Clifford algebra is a remarkably powerful tool to deal with the geometry of linear spaces \cite{Porteous_1995} but not a trivial one. Here we just review the properties of the neutral space $\R^{n,n}$ and of its Clifford algebra $\myClg{}{}{\R^{n,n}}$ at the heart of following results addressing the curious reader to quoted references.
%
%Clifford algebra is a remarkably powerful tool to deal with the geometry of linear spaces \cite{Porteous_1995} but not a trivial one. Here we just review the properties of the neutral space $\R^{n,n}$ and of its Clifford algebra $\myClg{}{}{\R^{n,n}}$ at the heart of following results addressing to quoted references for thorough treatment.
%
Clifford algebra is a remarkably powerful tool to deal with the geometry of linear spaces \cite{Porteous_1995} but not a trivial one. Here we just review the properties of the neutral space $\R^{n,n}$ and of its Clifford algebra $\myClg{}{}{\R^{n,n}}$ at the heart of following results. The interested reader has many choices to deepen the subject, here we mainly refer to the remarkable work of Ian Porteous \cite{Porteous_1995}.

%We review here the properties of the neutral space $\R^{n,n}$ and of its Clifford algebra $\myClg{}{}{\R^{n,n}}$ at the heart of following results. Clifford algebra is a remarkably powerful tool to deal with the geometry of linear spaces but not a trivial one and the interested reader has many choices to deepen this subject, here we refer to the remarkable work of Ian Porteous \cite{Porteous_1995}.

Clifford algebra contains a ``copy'' of the linear space so that it is customary to think at vectors as at algebra elements thus validating algebra product also for vectors. Clifford algebra admits a faithful and irreducible representation on any of its spinor spaces namely the linear spaces generated by its Minimal Left Ideals (MLI)%
\footnote{in an algebra $A$ a subset $\mySpinorS$ is a left ideal if for any $a \in A, \varphi \in \mySpinorS \implies a \varphi \in \mySpinorS$; it is minimal if it does not contain properly any other ideal. For example in matrix algebra the subset of matrices with only one nonzero column form a minimal left ideal.}%
. So the same element of the algebra may be thought both as a vector and as an endomorphism of spinor spaces.

In particular $\myClg{}{}{\R^{n,n}}$ is isomorphic to the algebra of real matrices $\R(2^n)$ \cite{Porteous_1995} spinors corresponding to columns. The $2 n$ generators of the algebra $\{ \mygen_{i} \}$ are the vectors of an orthonormal basis of the linear space $\R^{n,n}$. Within the algebra the familiar vector scalar product corresponds to an anticommutator%
\opt{margin_notes}{\mynote{mbh.note: put here bilinear forms (2) ?
}}%
\begin{equation}
\label{formula_generators}
%\mygen_i \mygen_j + \mygen_j \mygen_i := \anticomm{\mygen_i}{\mygen_j} = 2 \delta_{i j} (-1)^{i+1} \qquad i,j = 1,2, \ldots, 2 n
\mygen_i \mygen_j + \mygen_j \mygen_i := \anticomm{\mygen_i}{\mygen_j} = 2 \left\{ \begin{array}{l l}
\delta_{i j} & \mbox{for} \; i \le n \\
- \delta_{i j} & \mbox{for} \; i > n
\end{array} \right.
\qquad i,j = 1,2, \ldots, 2 n
\end{equation}
where $\delta_{i j}$ is the Kronecker delta: $1$ iff $i = j$ and $0$ otherwise. We define the Witt, or null, basis of $\R^{n,n}$:
\begin{equation}
\label{formula_Witt_basis}
\left\{ \begin{array}{l l l}
p_{i} & = & \frac{1}{2} \left( \mygen_{i} + \mygen_{i + n} \right) \\
q_{i} & = & \frac{1}{2} \left( \mygen_{i} - \mygen_{i + n} \right)
\end{array} \right.
%\left\{ \begin{array}{l l l}
%p_{i} & = & \frac{1}{2} \left( \mygen_{2i-1} + \mygen_{2i} \right) \\
%q_{i} & = & \frac{1}{2} \left( \mygen_{2i-1} - \mygen_{2i} \right)
%\end{array} \right.
%\Rightarrow
%\left\{\begin{array}{l l l}
%\mygen_{2i-1} & = & p_{i} + q_{i} \\
%\mygen_{2i} & = & p_{i} - q_{i}
%\end{array} \right.
\quad i = 1,2, \ldots, n
\end{equation}
that, with $\mygen_{i} \mygen_{j} = - \mygen_{j} \mygen_{i}$ for $i \ne j$, gives
\begin{equation}
\label{formula_Witt_basis_properties}
\anticomm{p_{i}}{p_{j}} = \anticomm{q_{i}}{q_{j}} = 0
\qquad
\anticomm{p_{i}}{q_{j}} = \delta_{i j}
\end{equation}
showing that all $p_i, q_i$ are mutually orthogonal, also to themselves, that implies $p_i^2 = q_i^2 = 0$ and are thus ``null'' vectors. Defining
\begin{equation}
\label{formula_Witt_decomposition}
\left\{ \begin{array}{l}
P = \my_span{p_1, p_2, \ldots, p_n} \\
Q = \my_span{q_1, q_2, \ldots, q_n}
\end{array} \right.
\end{equation}
it is easy to verify that any two vectors of $P$ (or $Q$) are null and orthogonal to each other and thus $P$ and $Q$ are two totally null subspaces of maximum dimension $n$, or Maximally Totally Null Planes (MTNP). $P$ and $Q$ constitute a Witt decomposition \cite{Porteous_1995} of $\R^{n,n}$ since $P \cap Q = \{ 0 \}$ and $P \oplus Q = \R^{n,n}$.

$\myClg{}{}{\R^{n,n}}$ is more easily manipulated exploiting the properties of its Extended Fock Basis (EFB, see \cite{Budinich_2016} and references therein) with which any algebra element is a linear superposition of spinors. The $2^{2 n}$ spinors forming EFB are given by all possible sequences
\begin{equation}
\label{EFB_def}
\psi = \psi_1 \psi_2 \cdots \psi_n \qquad \psi_i \in \{ q_i p_i, p_i q_i, p_i, q_i \} \qquad i = 1, 2, \ldots ,n \dotinformula
\end{equation}
Since $\mygen_{i} \mygen_{i + n} = q_i p_i - p_i q_i := \comm{q_i}{p_i}$ in EFB the identity $\Identity$ and the volume element $\omega$ (scalar and pseudoscalar) assume similar expressions \cite{Budinich_2016}:
\begin{equation}
\label{identity_omega}
\begin{array}{l l l}
\Identity & := & \anticomm{q_1}{p_1} \anticomm{q_2}{p_2} \cdots \anticomm{q_n}{p_n} \\
\omega & := & \mygen_1 \mygen_2 \cdots \mygen_{2 n} = \comm{q_1}{p_1} \comm{q_2}{p_2} \cdots \comm{q_n}{p_n} \dotinformula
\end{array}
\end{equation}

\noindent Any element of (\ref{EFB_def}) is a spinor and an element of one of the $2^n$ MLI $\mySpinorS$ of $\myClg{}{}{\R^{n,n}}$; all these MLI are equivalent%
\opt{margin_notes}{\mynote{mbh.ref: Porteous 1995 p. 133 B\&T p. ? %\cite[p. ?]{Benn_1987}
}}%
{} in the sense that each of them can carry a representation of the algebra; moreover the algebra, as a vector space, is the direct sum of these spinor spaces.
%In the isomorphic matrix algebra $\R(2^n)$ one spinor space is usually the linear space generated by one of the $2^n$ columns of the matrix.

Any spinor of one MLI $\mySpinorS$ is a linear combination of $2^n$ spinors (\ref{EFB_def}) that form the Fock basis $\myFockB$ of that spinor space \cite{BudinichP_1989, Budinich_2016}. Fock basis spinors are called simple and are in one to one correspondence with MTNP \cite{BudinichP_1989}. Given simple spinor $\psi \in \mySpinorS_s \subset \mySpinorS$ (\ref{EFB_def}) let $M(\psi)$ be its corresponding MTNP
\begin{equation}
\label{formula_MTNP_M_psi}
M(\psi) = \{ v \in \R^{n,n} : v \ne 0, v \psi = 0 \} \dotinformula
\end{equation}
For example given simple spinor $\psi = p_1 \; q_2 p_2 \; q_3 p_3$ of $\myClg{}{}{\R^{3,3}}$ for any $v \in \my_span{p_1, q_2, q_3}$ we easily get $v \psi = 0$ (\ref{formula_Witt_basis_properties}) and so $M(\psi) = \my_span{p_1, q_2, q_3}$. More in general in $\myClg{}{}{\R^{n,n}}$
\begin{equation}
\label{formula_MTNP_M_psi_2}
M(\psi) = \my_span{x_1, x_2, \ldots, x_n} \quad x_i = \left\{ \begin{array}{l l}
p_i & \mbox{iff} \; \psi_i = p_i, p_i q_i \\
q_i & \mbox{iff} \; \psi_i = q_i, q_i p_i
\end{array} \right.
\quad i = 1,2, \ldots, n
\end{equation}
and let $\mysetM{}$ be the set of the $2^n$ MTNP $M(\psi)$. Each MTNP of \mysetM{} is the span of $n$ null vectors obtained choosing one element from each of the $n$ couples $(p_i, q_i)$ (\ref{formula_Witt_basis}) \cite{BudinichP_1989, Budinich_2016}.

\mysetM{} is a subset of the larger set \mysnqG{} of all MTNP of $\R^{n,n}$, a semi-neutral quadric Grassmannian according to Porteous \cite{Porteous_1995}. \mysnqG{} in turn is isomorphic to the subgroup $\OO{n}$ of $\OO{n,n}$, moreover $\OO{n}$ acts transitively on \mysnqG{} and thus also on \mysetM{}.

We review these relations: seeing the neutral space $\R^{n,n}$ as $\R^n \times \R^n$ we can write its generic element as $(x,y)$ and then $(x,y)^2 = x^2 - y^2$ and so for any $x \in \R^n \times \{0\}$ and $t \in \OO{n}$ $(x,t(x))$ is a null vector since $(x,t(x))^2 = x^2 - t(x)^2 = 0$. All pairs $(x,t(x))$ form a MTNP of $\R^{n,n}$ that we indicate, with self-explanatory notation, as $(\Identity, t)$. Isometry $t \in \OO{n}$ establishes the quoted isomorphism since any MTNP of \mysnqG{} can be written as $(\Identity, t)$ \cite[Corollary~14.13]{Porteous_1995} and thus the one to one correspondence between MTNP of $\R^{n,n}$ and $t \in \OO{n}$ is formally established by the map
\opt{margin_notes}{\mynote{mbh.note: is this the Cayley chart ?}}%
\begin{equation}
\label{formula_bijection_Nn_O(n)}
\OO{n} \to \mysnqG{}; t \to (\Identity, t) \qquad \implies \qquad \mysnqG = \{ (\Identity, t) : t \in \OO{n} \} \dotinformula
\end{equation}

For example in this setting two generic null vectors of $P$ and $Q$ (\ref{formula_Witt_decomposition}) are respectively $(x, x)$ and $(x, -x)$ and in our notation we represent the whole MTNP $P$ and $Q$ (\ref{formula_Witt_decomposition}) with
\opt{margin_notes}{\mynote{mbh.ref: see pp. 55, 56}}%
\begin{equation}
\label{formula_P_Q_def}
\left\{ \begin{array}{l}
P = (\Identity, \Identity) \\
Q = (\Identity, -\Identity) \dotinformula
\end{array} \right.
\end{equation}
The action of $\OO{n}$ is transitive on \mysnqG{} since for any $t, u \in \OO{n}$, $(\Identity,ut) \in \mysnqG$ and the action of $\OO{n}$ is trivially transitive on $\OO{n}$.

\section{The SATisfiability problem in Clifford algebra}
\label{SAT_in_Cl(g)}
We summarize \SAT{} problems \cite{Knuth_2015} formulation in $\myClg{}{}{\R^{n,n}}$ presented more extensively in other papers \cite{Budinich_2017, Budinich_2021} to which we address the interested reader.

In a nutshell: the conjunctive normal form of a $k\SAT{}$ problem $\myBooleanS$ with $n$ \myBooleanv{}s $\literal_1, \literal_2, \cdots, \literal_n$ and $m$ clauses ${\cal C}_j \equiv (\literal_{j_1} \lor \literal_{j_2} \lor \cdots \lor \literal_{j_k})$ (we use $\equiv$ for logical equivalence to avoid confusion with algebraic equality $=$) is
$$
\myBooleanS \equiv {\cal C}_1 \land {\cal C}_2 \land \cdots \land {\cal C}_m \dotinformula
$$
Formally a solution is either an assignment for the $n$ \myBooleanv{}s $\literal_i \in \{ \myBooleanT, \myBooleanF \}$ (true, false) that make $\myBooleanS \equiv \myBooleanT$ or a proof that none exists. We can formulate $\myBooleanS$ in $\myClg{}{}{\R^{n,n}}$ with the following substitutions that replace Boolean expressions with algebraic ones in $\myClg{}{}{\R^{n,n}}$ ($\myconjugate{\literal}_i$ stands for $\lnot \literal_i$)
\begin{equation}
\label{Boolean_substitutions}
\begin{array}{lll}
\myBooleanF & \to & 0 \\
\myBooleanT & \to & \Identity \\
\literal_i & \to & q_i p_i \\
%L_1 \land L_2 \to l_1 l_2 \\
%\lnot \literal_i \equiv
\myconjugate{\literal}_i & \to & \Identity - q_i p_i = p_i q_i\\
\literal_i \land \literal_j & \to & q_i p_i \; q_j p_j
\end{array}
\end{equation}
$p_i$ and $q_i$ being vectors of the Witt basis (\ref{formula_Witt_basis}). With (\ref{formula_Witt_basis_properties}) we easily get
$$
q_i p_i \; q_i p_i = q_i p_i \quad p_i q_i \; p_i q_i = p_i q_i \quad q_i p_i \; p_i q_i = p_i q_i \; q_i p_i = 0 \quad q_i p_i \; q_j p_j = q_j p_j \; q_i p_i
$$
that shows that $q_i p_i$ and $p_i q_i$ are part of a family of orthogonal, commuting, idempotents and with a simple exercise we see that all $\myClg{}{}{\R^{n,n}}$ elements in (\ref{Boolean_substitutions}) are idempotents.
% and, with the associations:
%$$
%\literal_i \to q_i p_i \qquad \myconjugate{\literal}_i \to p_i q_i \qquad \myBooleanF \to 0
%$$
%and mapping the logical AND to Clifford product, we recognize in previous relations the Boolean relations
%$$
%\literal_i \land \myconjugate{\literal}_i \equiv \myconjugate{\literal}_i \land \literal_i \equiv \myBooleanF \quad \literal_i \land \literal_i \equiv \literal_i \quad \myconjugate{\literal}_i \land \myconjugate{\literal}_i \equiv \myconjugate{\literal}_i \quad \literal_i \land \literal_j \equiv \literal_j \land \literal_i \dotinformula
%$$
%From now on we will use $\literal_i$ and $\myconjugate{\literal}_i$ also in $\myClg{}{}{\R^{n,n}}$ meaning respectively $q_i p_i$ and $p_i q_i$ and Clifford product will stand for logical AND.
%
%To lift these promising relations to a full Boolean algebra we must add further structure.
%
For example $q_i p_i \; q_i p_i = q_i p_i$ stands for the logical relation $\literal_i \land \literal_i \equiv \literal_i$. From now on we will use $\literal_i$ and $\myconjugate{\literal}_i$ also in $\myClg{}{}{\R^{n,n}}$ meaning respectively $q_i p_i$ and $p_i q_i$ and Clifford product will stand for logical AND $\land$. In this setting we can prove \cite{Budinich_2017} the following result.
\begin{MS_Proposition}
\label{SAT_in_Cl_2}
Given a SAT problem $\myBooleanS$ with $m$ clauses ${\cal C}_j \equiv (\literal_{j_1} \lor \literal_{j_2} \lor \cdots \lor \literal_{j_k})$, for each clause let $z_j := \myconjugate{\literal}_{j_1} \myconjugate{\literal}_{j_2} \cdots \myconjugate{\literal}_{j_k}$, then $\myBooleanS$ is unsatisfiable if and only if, for the corresponding algebraic expression of $\myClg{}{}{\R^{n,n}}$
\begin{equation}
\label{formula_SAT_EFB_2}
S = \prod_{j = 1}^m (\Identity - z_j) = 0 \dotinformula
\end{equation}
\end{MS_Proposition}

We remark that $z_j$ represents the unique assignment (more on this below) of the $k$ literals of ${\cal C}_j$ that give ${\cal C}_j \equiv \myBooleanF$ and thus $\Identity - z_j$ substantially means all possible assignments of the \myBooleanv{}s \emph{but} $z_j$ that allows to grasp intuitively the rationale behind the algebraic expression of $S$ (\ref{formula_SAT_EFB_2}). From now on we will represent clauses ${\cal C}_j$ only in form $z_j$.

A $1\SAT$ problem is just a logical AND of $m$ literals. For both assignments of $\literal_i$, $\literal_i \land \myconjugate{\literal}_i \equiv \myBooleanF$ and thus the presence of a literal together with its logical complement is a necessary and sufficient condition for making a $1\SAT$ formula unsatisfiable. We can interpret a satisfiable $1\SAT$ formula as an assignment for its variables since there is only one assignment of its variables that makes it $\myBooleanT$ and that can be read scanning the formula; in the sequel we will freely use $1\SAT$ formulas for assignments of variables and we will switch between the two forms as and when it suits us. We also easily see that $1\SAT$ formulas are idempotents.

Since $\myClg{}{}{\R^{n,n}}$ is a simple algebra, the unit element of the algebra is the sum of $2^n$ primitive (indecomposable) idempotents $\myprimidemp_i$
\begin{equation}
\label{formula_identity_def}
\Identity = \sum_{i = 1}^{2^n} \myprimidemp_i = \prod_{j = 1}^{n} \anticomm{q_j}{p_j}
\end{equation}
where the product of $n$ anticommutators is its expression in EFB (\ref{identity_omega}). The full expansion of these anticommutators contains $2^n$ terms each term being one of the primitive idempotents \emph{and} a simple spinor (\ref{EFB_def}). At this point it is manifest that the $2^n$ primitive idempotents $\myprimidemp_i$ of the expansion (\ref{formula_identity_def}) are in one to one correspondence with the possible $2^n$ $1\SAT$ formulas of the $n$ literals (Boolean atoms), for example:
$$
\literal_1 \myconjugate{\literal}_2 \cdots \literal_n \to q_1 p_1 \; p_2 q_2 \cdots q_n p_n \dotinformula
$$
Given an assignment of its $n$ \myBooleanv{}s, \eg $\literal_1 \myconjugate{\literal}_2 \cdots \literal_n$, $\myBooleanS \equiv \myBooleanF$ if and only if $\literal_1 \myconjugate{\literal}_2 \cdots \literal_n \land \myBooleanS \equiv \myBooleanF$, becoming $\literal_1 \myconjugate{\literal}_2 \cdots \literal_n \, S= 0$ in Clifford algebra. We interpret these formulas as the \emph{substitution} of the only assignment satisfying $\literal_1 \myconjugate{\literal}_2 \cdots \literal_n$ into $\myBooleanS$. By (\ref{formula_SAT_EFB_2}) $\literal_1 \myconjugate{\literal}_2 \cdots \literal_n \, S = 0$ if and only if there exists a clause $z_j$ such that $\literal_1 \myconjugate{\literal}_2 \cdots \literal_n (\Identity - z_j) = 0$ namely $\literal_1 \myconjugate{\literal}_2 \cdots \literal_n = \literal_1 \myconjugate{\literal}_2 \cdots \literal_n z_j$ \cite{Budinich_2017}.

It is instructive to derive the full expression of \eg $q_1 p_1$ directly from EFB formalism; with (\ref{formula_identity_def})
\begin{equation}
\label{literal_projection}
q_1 p_1 = q_1 p_1 \Identity = q_1 p_1 \prod_{j = 2}^{n} \anticomm{q_j}{p_j}
\end{equation}
since $q_1 p_1 \anticomm{q_1}{p_1} = q_1 p_1$ and the full expansion is a sum of $2^{n - 1}$ EFB terms that are all primitive idempotents and thus $q_1 p_1$ is an idempotent the sum being precisely the expansion as a sum of the primitive idempotents $\myprimidemp_i$. From the logical viewpoint this can be interpreted as the property that given the $1\SAT$ formula $\literal_1$ the other, unspecified, $n-1$ \myBooleanv{}s $\literal_2, \ldots, \literal_n$ are free to take all possible $2^{n - 1}$ values or, more technically, that $\literal_1$ has a \emph{full} disjunctive normal form made of $2^{n - 1}$ Boolean atoms.

\section{\SAT{} Clauses and Totally Null Planes of $\R^{n,n}$}
\label{clauses_TNP}
We are now ready to exploit the \SAT{} formulation in Clifford algebra to transform a \SAT{} problem $S$ into a geometric problem of null subspaces of $\R^{n,n}$. By previous findings any assignment $\literal_1 \literal_2 \cdots \literal_n$ can represent:
\begin{itemize}
\item a 1\SAT{} formula and a Boolean atom,
\item an element of a MLI of $\myClg{}{}{\R^{n,n}}$,
\item a simple spinor in $\mySpinorS_s$ and thus, by (\ref{formula_MTNP_M_psi_2}),
\item a MTNP $M(\literal_1 \literal_2 \cdots \literal_n) \in \mysetM{}$%
\footnote{A technical remark about this passage: strictly any assignment $\literal_1 \literal_2 \cdots \literal_n$ belongs to a different MLI of the $2^n$ of $\myClg{}{}{\R^{n,n}}$ but it may be ``projected'' to any other MLI while its associated MTNP remains the same \cite{Budinich_2016}.}%
.
\end{itemize}

Along the same path also a clause in the form $z_j$ of Proposition~\ref{SAT_in_Cl_2} defines a Totally Null Plane (TNP) that in general is not maximal having dimension less than $n$ given that $z_j$ is an assignment of $k < n$ \myBooleanv{}s. Unfortunately whereas the correspondence between simple spinors and MTNP of \mysetM{} is one to one the case of generic spinors and their associated TNP is more intricate. In this case we can prove the following.
\opt{margin_notes}{\mynote{mbh.note: I should perhaps prove that the map $z_j \to {\cal N}_k$ is injective}}%
\begin{MS_Proposition}
\label{M_z_j}
Any clause $z_j$ determines uniquely
\opt{margin_notes}{\mynote{mbh.note: is unicity proved ?}}%
a spinor $\psi_{z_j}$, in general not simple. If $z_j$ has $k < n-2$ literals the spinor $\psi_{z_j}$ induces $k$TNP $M(z_j)$ of dimension $k$ and given by (\ref{formula_MTNP_M_psi_2}) applied to the $k$ literals forming $z_j$.
%
%\bigskip
%
%
%
%
%Given a clause $z_j$ made by $k < n-2$ literals it determines uniquely
%%
%\opt{margin_notes}{\mynote{mbh.note: is unicity proved ?}}%
%%
%a spinor $\psi_{z_j}$ which TNP $M(\psi_{z_j}) = M(z_j)$ has dimension $k$ and is the span of the null vectors appearing at the first place in the literals that form $z_j$.
\end{MS_Proposition}

\begin{proof}
Applying to $z_j$ the same procedure used in expansion (\ref{literal_projection}) we get
\begin{equation}
\label{formula_z_j_EFB}
z_j = z_j \Identity= \myconjugate{\literal}_{j_1} \myconjugate{\literal}_{j_2} \cdots \myconjugate{\literal}_{j_k} \Identity = \myconjugate{\literal}_{j_1} \myconjugate{\literal}_{j_2} \cdots \myconjugate{\literal}_{j_k} \prod_{i \notin j_k} \anticomm{q_i}{p_i}
\end{equation}
that expands in a sum of $2^{n-k}$ elements of $\myClg{}{}{\R^{n,n}}$. Each of these $2^{n-k}$ elements belongs to a different MLI \cite{Budinich_2016} but they can all be ``projected'' in a unique MLI where they form the Fock basis expansion of a spinor $\psi_{z_j}$
\begin{equation}
\label{formula_psi_z_j_EFB}
\psi_{z_j} = \psi_{j_1} \psi _{j_2} \cdots \psi_{j_k} \Psi
\end{equation}
where $\psi_{j_1} \psi _{j_2} \cdots \psi_{j_k}$ and $\Psi$ are respectively the projections of $\myconjugate{\literal}_{j_1} \myconjugate{\literal}_{j_2} \cdots \myconjugate{\literal}_{j_k}$ and of the totally symmetric product of anticommutators $\prod_{i \notin j_k} \anticomm{q_i}{p_i}$ in one MLI, namely spinor space $\mySpinorS$. By Theorem~1 of \cite{Budinich_2014} if $n - k > 2$ then spinor $\psi_{z_j}$ has an associated $k$TNP $M(z_j)$ formed by all and only the $k$ null vectors given by (\ref{formula_MTNP_M_psi_2}) that annihilate its first part $\psi_{j_1} \psi _{j_2} \cdots \psi_{j_k}$ given that there are no vectors in $\R^{n,n}$ that annihilate $\Psi$ (see \cite{Budinich_2014} for details).
\end{proof}
It is thus appropriate, for $k < n-2$, this usually being the case in ``real life'' \SAT{} problems, to associate also to a clause $z_j$ its $k$TNP $M(z_j)$.%
\opt{margin_notes}{\mynote{mbh.note: for $k \ge n-2$ are the following Propositions affected ?}}%
{} More precisely if $r$ of the $k$ literals of $z_j$ appear in plain form and $k-r$ in complementary form, then with (\ref{Boolean_substitutions}) and (\ref{formula_MTNP_M_psi_2}) we get
\begin{equation}
\label{M_z_j_def}
M(z_j) = \my_span{q_{j_1}, \ldots, q_{j_r}, p_{j_1}, \ldots, p_{j_{k-r}}} \dotinformula
\end{equation}
Moreover each of the $2^{n-k}$ elements of the expansion of $\psi_{z_j}$ (\ref{formula_psi_z_j_EFB}) is a simple spinor which MTNP contains the $k$ null vectors of $M(z_j)$ with a combination of the other possible $n-k$ null vectors; in this fashion $M(z_j)$ is contained in all the $2^{n-k}$ MTNP of this expansion and is thus the $k$TNP corresponding precisely to the intersection of these $2^{n-k}$ MTNP of \mysetM{}; we have thus proved the following Corollary.
\begin{MS_Corollary}
\label{M_z_j_coro}
The $k$TNP $M(z_j)$ (\ref{M_z_j_def}), associated to a clause $z_j$ with $k < n-2$ literals, is the intersection of the $2^{n-k}$ MTNP of \mysetM{} associated to the Fock basis expansion of spinor $\psi_{z_j}$ (\ref{formula_psi_z_j_EFB}).
%
%
%\bigskip
%
%Given a clause $z_j$ made by $k < n-2$ literals then its $k$TNP $M(z_j)$ (\ref{M_z_j_def}) is the intersection of the $2^{n-k}$ MTNP of \mysetM{} forming $\psi_{z_j}$ (\ref{formula_psi_z_j_EFB}).
\end{MS_Corollary}

We define a clause $z_j$ and an assignment $\literal_1 \literal_2 \cdots \literal_n$ to be \emph{compatible} if, equivalently,
\opt{margin_notes}{\mynote{mbh.ref: see log p. 731.1}}%
\begin{equation}
\label{formula_compatible_def}
\left\{ \begin{array}{l}
z_j \quad \mbox{literals appear in the same form also in} \quad \literal_1 \literal_2 \cdots \literal_n \\
\literal_1 \literal_2 \cdots \literal_n z_j = \literal_1 \literal_2 \cdots \literal_n \;\; \iff \;\; \literal_1 \literal_2 \cdots \literal_n (\Identity - z_j) = 0 \\
M(z_j) \subseteq M(\literal_1 \literal_2 \cdots \literal_n) \; \iff \; M(z_j) \cap M(\literal_1 \literal_2 \cdots \literal_n) = M(z_j)
\end{array} \right.
\end{equation}
and with Corollary~\ref{M_z_j_coro} it is simple to verify that each definition implies the successive circularly.

\bigskip

We continue showing that if all $2^n$ assignments $\literal_1 \literal_2 \cdots \literal_n$ are compatible with at least one clause $z_j$ then a \SAT{} problem is unsatisfiable.
\begin{MS_Proposition}
\label{SAT_in_TNP}
A \SAT{} problem $\myBooleanS$ is unsatisfiable if and only if, for any of the $2^n$ assignments $\literal_1 \literal_2 \cdots \literal_n$, at least one of the $m$ clauses $z_j$ of $S$ is compatible (\ref{formula_compatible_def}) with the assignment.
\end{MS_Proposition}

\begin{proof}
Let $\myBooleanS$ be unsatisfiable: for any assignment $\literal_1 \literal_2 \cdots \literal_n$ then $\literal_1 \literal_2 \cdots \literal_n S = 0$ as shown in the second part of Section~\ref{SAT_in_Cl(g)}. Given $S$ expression (\ref{formula_SAT_EFB_2}) this happens if and only if there exists at least one $z_j$ such that $\literal_1 \literal_2 \cdots \literal_n (\Identity - z_j) = 0$, namely iff $z_j$ is compatible with $\literal_1 \literal_2 \cdots \literal_n$.

Conversely let $S$ be such that for any $\literal_1 \literal_2 \cdots \literal_n$ there exists at least one clause $z_j$ such that $\literal_1 \literal_2 \cdots \literal_n (\Identity - z_j) = 0$; it follows that for any assignment $\literal_1 \literal_2 \cdots \literal_n S = 0$ that proves that $S$ is unsatisfiable.
\end{proof}

We remark that the complementary statement that given an assignment $\literal_1 \literal_2 \cdots \literal_n$ there are clauses such that $\cup_j M(z_j) = M(\literal_1 \literal_2 \cdots \literal_n)$ is not true; for example $S = \literal_1 \myconjugate{\literal}_1$ is unsatisfiable but $M(\myconjugate{\literal}_1) \cup M(\literal_1)$ is not even a TNP.

With this result a \SAT{} problem is transformed into a geometrical problem of MTNP in $\R^{n,n}$ and in what follows we push forward in this direction.

\section{\SAT{} Clauses and discrete isometries of $\OO{n}$}
\label{clauses_isometries}
We have just seen that a clause $z_j$ of a \SAT{} problem $\myBooleanS$ defines a $k$TNP $M(z_j)$ but since \mysnqG{} and $\OO{n}$ are isomorphic (\ref{formula_bijection_Nn_O(n)}) it is not surprising that a clause may induce also an isometry $t_j \in \OO{n}$. We begin showing that $\myBooleanS$ is unsatisfiable if and only if the isometries $t_j$ induced by its clauses $z_j$ form a cover of an abelian subgroup of $\OO{n}$.

We start investigating bijection (\ref{formula_bijection_Nn_O(n)}) when restricted to the subset $\mysetM \subset \mysnqG{}$ and we take $P = (\Identity, \Identity)$ (\ref{formula_P_Q_def}) as our ``reference'' MTNP of \mysetM{}.
\opt{margin_notes}{\mynote{mbh.note: in physics language this corresponds to the choice of the ``vacuum'' spinor}}%

Let $t_i$ be the isometry (actually an involution) that inverts the timelike vector $\mygen_{i + n}$, namely
\begin{equation}
\label{t_i_def}
t_i(\mygen_{j}) =
\left\{ \begin{array}{l l l}
-\mygen_{j} & \quad \mbox{for} \quad j = i + n\\
\mygen_{j} & \quad \mbox{otherwise}
\end{array} \right.
\quad i = 1,2, \ldots, n \quad j = 1,2, \ldots, 2n
\end{equation}
its action on the Witt basis (\ref{formula_Witt_basis}) is to exchange the null vectors $p_i$ and $q_i$. It follows that the inversion of a certain subset of timelike vectors exchanges the corresponding null vectors $p_i$ with $q_i$ and vice-versa. It is thus clear that starting from $P$ we can obtain any element of \mysetM{} by the corresponding inversion of a subset of the $n$ timelike vectors $\mygen_{i + n}$. Each of these isometries acts on the (timelike) subspace $\{0\} \times \R^n$ of $\R^{n,n}$ and, in the matrix representation of $\OO{n}$, is a diagonal matrix $\lambda \in \R(n)$ with $\pm 1$ on the diagonal and all these matrices form the group $\mbox{O}(1) \times \mbox{O}(1) \cdots \times \mbox{O}(1) = \stackrel{n}{\times} \mbox{O}(1) := \O1n$ that is immediate to get remembering that $\mbox{O}(1) = \{ \pm 1 \}$.%
\opt{margin_notes}{\mynote{mbh.ref: remember that the tensor product is associative Proposition~11.5 of {Porteous 1995}; $\O1n$ not normal $t \lambda t^{-1} \notin \O1n$}}%
{} \O1n is a discrete, abelian, subgroup of $\OO{n}$ and its elements are involutions since $\lambda^2 = \Identity$. It is thus clear that for any assignment $\literal_1 \literal_2 \cdots \literal_n$ there exists a unique involution $\lambda$ giving
\begin{equation}
\label{formula_assignment_involution}
M(\literal_1 \literal_2 \cdots \literal_n) = (\Identity, \lambda) \qquad \lambda \in \O1n
\end{equation}
and $(\Identity, \lambda)$ can also be seen as obtained by the action of $\lambda$ on $P$ since $\lambda \Identity = \lambda$
so we proved constructively:
\begin{MS_Proposition}
\label{isomorphism_restricted}
The isomorphism (\ref{formula_bijection_Nn_O(n)}) when restricted to the subgroup \O1n of $\OO{n}$ has for image $\mysetM \subset \mysnqG$ and \O1n acts transitively on \mysetM{}.
\end{MS_Proposition}

We can thus enrich the list of possible interpretations of an assignment $\literal_1 \literal_2 \cdots \literal_n$ made in Section~\ref{clauses_TNP} with that of an involution $\lambda \in \O1n$ (\ref{formula_assignment_involution}). We are now ready to define the isometry $\lambda_j \in \O1n$ associated to a clause $z_j = \literal_{j_1} \literal_{j_2} \cdots \literal_{j_k}$ with $k$ literals;
%let $P_j = \my_span{p_{j_1}, p_{j_2}, \ldots, p_{j_k}}$ be the corresponding $k$TNP subspace of $P$
let $r$ of these $k$ literals appear in plain form and $k-r$ in complementary form, then $M(z_j)$ is as in (\ref{M_z_j_def}). Using (\ref{t_i_def}) let
\begin{equation}
\label{formula_Lambda_j_def}
\lambda_j = t_{j_1} \cdots t_{j_r}
\end{equation}
with the understanding that $\lambda_j = \Identity$ if $r = 0$, then it is simple to verify that %$M(z_j) = \lambda_j (P_j)$ and thus
\begin{equation}
\label{formula_Lambda_j_action}
M(z_j) \subseteq (\Identity, \lambda_j)
\end{equation}
namely that the isometry $\lambda_j$ (\ref{formula_Lambda_j_def}) is a MTNP of \mysetM{} that has $M(z_j)$ as a subspace: in practice $\lambda_j$ exchanges $r$ of the $n$ $p_i$ of $P$ with $q_i$. The definition of $\lambda_j$ (\ref{formula_Lambda_j_def}) satisfying (\ref{formula_Lambda_j_action}) is not unique since we can freely add to $\lambda_j$ (\ref{formula_Lambda_j_def}) any subset of the $n-k$ involutions $t_i$ whose indexes do not appear in $z_j$ and (\ref{formula_Lambda_j_action}) continues to hold. There are thus $2^{n-k}$ different $\lambda_j$ all satisfying (\ref{formula_Lambda_j_action}) and they are the $2^{n-k}$ MTNP $M(z_j)$ of Corollary~\ref{M_z_j_coro} corresponding to the EFB expansion (\ref{formula_psi_z_j_EFB}). We just remark that all $\lambda_j$ are involutions of \O1n.

We define the set of involutions induced by a clause $z_j$ as
\begin{equation}
\label{formula_cal_T_j_def}
{\cal T}_j' := \{\lambda_j \in \O1n : M(z_j) \subseteq (\Identity, \lambda_j) \}
\end{equation}
namely the set of $2^{n-k}$ elements of \O1n that satisfy (\ref{formula_Lambda_j_action}) and also, by (\ref{formula_assignment_involution}) and (\ref{formula_compatible_def}), the set of $2^{n-k}$ assignments compatible with $z_j$.
\opt{margin_notes}{\mynote{mbh.note: and also a subset of the Fock basis}}%
\begin{MS_Proposition}
\label{SAT_in_n_O(1)}
A given \SAT{} problem $\myBooleanS$ with m clauses $z_j$ is unsatisfiable if and only if the isometries induced by its clauses (\ref{formula_cal_T_j_def}) form a cover for \O1n:%
\opt{margin_notes}{\mynote{mbh.ref: cover defined at p.~613 of logbook. Would it be more precise to use ``atlas'' ?}}%
\begin{equation}
\label{formula_SAT_in_n_O(1)}
\cup_{j = 1}^m {\cal T}_j' = \O1n \dotinformula
\end{equation}
\end{MS_Proposition}
\begin{proof}
Let $\myBooleanS$ be unsatisfiable: for any assignment $\lambda$ (\ref{formula_assignment_involution}) by Proposition~\ref{SAT_in_TNP} there exists at least one compatible clause $z_j$ (\ref{formula_compatible_def}) and thus $\lambda \in {\cal T}_j'$. Conversely from (\ref{formula_SAT_in_n_O(1)}) any $\lambda \in \O1n$ belongs to at least one ${\cal T}_j'$ of clause $z_j$ and thus assignment $\lambda$ (\ref{formula_assignment_involution}) is compatible with $z_j$ and $\myBooleanS$ is unsatisfiable.
\end{proof}

\section{\SAT{} Clauses and continuous isometries of $\OO{n}$}
\label{continuous_isometries}
In the last step we show that when a \SAT{} problem $\myBooleanS$ is unsatisfiable the isometries induced by its clauses not only form a cover of \O1n (\ref{formula_SAT_in_n_O(1)}) but also of its parent group $\OO{n}$.

We start extending the definition of isometries induced by a clause (\ref{formula_cal_T_j_def}) to
\begin{equation}
\label{formula_cal_T_j_def2}
{\cal T}_j := \{t \in \OO{n} : M(z_j) \subseteq (\Identity, t)\} \subset \OO{n}
\end{equation}
this being an obvious generalization of (\ref{formula_cal_T_j_def}), moreover $\O1n \subset \OO{n}$ implies ${\cal T}_j' \subseteq {\cal T}_j$.
To proceed we need:
\begin{MS_Proposition}
\opt{margin_notes}{\mynote{mbh.note: log pp. 774, 783. Needs definitions: $\{\psi_\lambda \}$, $\psi = \sum_\lambda \alpha_\lambda \psi_\lambda$ etc. Commented there is here another proof that exploits Fock basis, too quickly defined in Section~\ref{neutral_R_nn}.}}%
\label{prop_Tau_j_equivalence}
Given any clause $z_j$ and ${\cal T}_j'$ (\ref{formula_cal_T_j_def}) the three following definitions of ${\cal T}_j \subset \OO{n}$ are equivalent
$$
{\cal T}_j =
\left\{ \begin{array}{l}
\{ t \in \OO{n} : M(z_j) \subseteq (\Identity, t) \} \subset \OO{n} \\
\{ \psi_t \in \my_span{\psi_\lambda} : \lambda \in {\cal T}_j' \subset \myFockB \mbox{ and } \psi_t \in \mySpinorS_s \mbox{ with } M(\psi_t) = (\Identity, t) \} \subset \mySpinorS_s \\
\{ t \in \OO{n} : (\Identity, t) = M(\psi_t) \mbox{ for } \psi_t \mbox{ as above} \}
\end{array} \right.
$$
\end{MS_Proposition}

%another proof
\begin{proof}
Assuming (\ref{formula_cal_T_j_def2}) let \eg $M(z_j) = \my_span{p_r, q_s, p_u}$ (\ref{M_z_j_def}) and thus for any $t \in {\cal T}_j$, $\my_span{p_r, q_s, p_u} \subseteq (\Identity, t)$ and for simple spinor $\psi_t$ such that $M(\psi_t) = (\Identity, t)$, $\my_span{p_r, q_s, p_u} \subseteq M(\psi_t)$. Also for any $\lambda_j \in {\cal T}_j'$, $\my_span{p_r, q_s, p_u} \subseteq (\Identity, \lambda_j)$ and all these spinors can be expanded in Fock basis $\myFockB$ of $\mySpinorS$ and any of these $\psi_\lambda \in \myFockB$ are identified by null vectors $p_r, q_s, p_u$ that are exactly those having $\lambda_j \in {\cal T}_j'$.
\opt{margin_notes}{\mynote{mbh.note: last part is confused, do not prove = of different defs, assumes Fock basis formalism}}%
\end{proof}
%
%
%
%\begin{proof}
%We show that each definition implies the next one, circularly. Given first definition the second one is the simple transposition in $\OO{n}$ of the bijection between simple spinors and $\OO{n}$ \cite{BudinichP_1989}. Given second definition we know that for any simple spinor $\psi$ than for any $v \in M(z_j)$ $v \psi = 0$ and so $v \in M(\psi)$ that implies $M(z_j) \subset M(\psi) = (\Identity, t)$ thus proving that the second definition implies the third one. Assuming the third definition we assume e.g. $M(z_j) = \my_span{p_r, q_s, p_u}$ and thus for any $t \in {\cal T}_j$ we know that $\my_span{p_r, q_s, p_u} \subseteq (\Identity, t)$ and thus for any $\psi \in \mySpinorS_s$ such that $M(\psi) = (\Identity, t)$ we have $\my_span{p_r, q_s, p_u} \subseteq M(\psi)$ and all these spinors can be expanded in Fock basis $\myFockB$ of $\mySpinorS$ and any of these $\psi_\lambda \in \myFockB$ are identified by null vectors $p_r, q_s, p_u$ that are exactly those having $\lambda \in {\cal T}_j'$.
%\end{proof}
%
\begin{MS_lemma}
\opt{margin_notes}{\mynote{mbh.note: log pp. 774, 783, commented there is another proof that uses Fock basis and thus needs the sum of ${\cal T}_j$}}%
\label{lemma_T_T'_intersection}
For any ${\cal T}_j$ then ${\cal T}_j \cap \O1n = {\cal T}_j'$\/.
\end{MS_lemma}
\begin{proof}
We already saw that ${\cal T}_j' \subseteq \O1n$ and ${\cal T}_j' \subseteq {\cal T}_j$ so we need just to prove that for any $\lambda' \in \O1n$, $\lambda' \in {\cal T}_j$ implies $\lambda' \in {\cal T}_j'$. Let us suppose the contrary: by Corollary~\ref{M_z_j_coro} there is at least one \eg $q_{j_x}$ such that $q_{j_x} \in M(z_j)$ while $p_{j_x} \in (\Identity, \lambda')$ (or viceversa) (\ref{M_z_j_def}). So we would have $\lambda' \in {\cal T}_j$ but $M(z_j) \not\subset (\Identity, \lambda')$ that is a contraddiction (\ref{formula_cal_T_j_def2}).
\end{proof}
%
%\begin{proof}
%Given bijection between $\O1n$ and Fock basis $\myFockB$ (\ref{formula_assignment_involution}) \cite[WRONG}{Budinich_2020} then for any $\psi_\lambda \in F : \lambda \in {\cal T}_j'$ is also in ${\cal T}_j$. Should any $\psi_\omega \in F : \omega \notin {\cal T}_j'$ be in ${\cal T}_j$ by ${\cal T}_j$ definition $\psi_\omega = \sum_\lambda \alpha_\lambda \psi_\lambda$ for $\lambda \in {\cal T}_j'$ but this is impossible since $F$ is a basis of $S$ with linearly independent elements \cite{BudinichP_1989, Budinich_2016}.
%\end{proof}

Any MTNP of \mysnqG{} is a subspace of $\R^{n,n}$ and two of them necessarily have an intersection of dimension $r$ with $0 \le r \le n$, their \emph{incidence}, and to proceed we need a crucial property of simple spinors \cite{BudinichP_1989} that we reproduce here, slightly adapted to our needs:
\begin{MS_Proposition}
\label{prop_5_BudinichP_1989}
Given any two linearly independent simple spinors $\psi, \phi \in \mySpinorS_s$ then their linear combinations $\alpha \psi + \beta \phi$ ($\alpha, \beta \in \R$) are simple if and only if the incidence of their associated MTNP is $n-2$ namely
$$
\dim M(\psi) \cap M(\phi) = n-2
$$
and then $M(\psi) \cap M(\psi + \phi) = M(\phi) \cap M(\psi + \phi) = M(\psi) \cap M(\phi)$.
\end{MS_Proposition}

We can thus define the set of simple spinors linear combinations of spinors coming from different subsets ${\cal T}_j, {\cal T}_k$ induced by clauses (\ref{formula_cal_T_j_def2})
\begin{equation}
\label{formula_Tau_addition}
{\cal T}_j + {\cal T}_k := \{ \psi = \alpha \psi_j + \beta \psi_k : \psi \in \mySpinorS_s; \psi_j \in {\cal T}_j; \psi_k \in {\cal T}_k; \alpha, \beta \in \R \}
\end{equation}
for example in the $k = n = 2$ \SAT{} problem given by the two clauses $z_1 = \literal_1 \literal_2$, $z_2 = \myconjugate{\literal}_1 \myconjugate{\literal}_2$ we easily get ${\cal T}_1 = {\cal T}_1' = \{ \Identity_2 \} \subset \SO{2}$ and ${\cal T}_2 = {\cal T}_2' = \{ - \Identity_2 \} \subset \SO{2}$ and ${\cal T}_1 + {\cal T}_2 = \SO{2}$ since $\psi_1 = p_1 q_1 p_2 q_2, \psi_2 = q_1 q_2$ with $M(\psi_1) = \my_span{p_1, p_2}$ and $M(\psi_2) = \my_span{q_1, q_2}$ of incidence $n - 2 = 0$ that satisfy Proposition~\ref{prop_5_BudinichP_1989} and we easily get
\opt{margin_notes}{\mynote{mbh.note: easily ? at least a ref; also $\theta/2$ ? see pp. 759, 777}}%
that for any $\theta \in [0, 2 \pi)$, $\psi = \cos \frac{\theta}{2} \psi_1 + \sin \frac{\theta}{2} \psi_2 \in \mySpinorS_s$ and that $M(\psi) = (\Identity, t_\theta)$ with $t_\theta = \left(\begin{array}{r r} \cos \theta & - \sin \theta \\ \sin \theta & \cos \theta \end{array}\right)$, namely $\SO{2}$.

We remark that neither ${\cal T}_j$ nor ${\cal T}_j + {\cal T}_k$ are linear subspaces of $\OO{n}$ since for any $t_1, t_2 \in {\cal T}_j$, with incidence of their MTNP different from $n-2$, their linear combination is not a MTNP and thus can't be in ${\cal T}_j$.
\opt{margin_notes}{\mynote{mbh.note: log p.806}}%
Moreover the sum operation of (\ref{formula_Tau_addition}) does \emph{not} define an operation between sets ${\cal T}$'s and with $\sum_j {\cal T}_j$ notation we just indicate the set of spinors with elements taken from sets ${\cal T}$'s and with this caveat we can generalize the simple $\SO{2}$ example to arbitrary $n$.

\begin{MS_theorem}
\opt{margin_notes}{\mynote{mbh.note: log pp. 773', 787, 795}}%
\label{SAT_in_O(n)}
A given \SAT{} problem $\myBooleanS$ with $n$ \myBooleanv{}s is unsatisfiable if and only if the isometries induced by its $m$ clauses (\ref{formula_cal_T_j_def2}) form a cover for $\OO{n}$:
\begin{equation}
\label{formula_SAT_in_O(n)}
\sum_{j = 1}^m {\cal T}_j = \OO{n} \dotinformula
\end{equation}
\end{MS_theorem}
\begin{proof}
Let $\myBooleanS$ be unsatisfiable, by Proposition~\ref{SAT_in_n_O(1)} $\cup_{j = 1}^m {\cal T}_j' = \O1n$ namely $\cup_{j = 1}^m {\cal T}_j' = \myFockB$
\opt{margin_notes}{\mynote{mbh.note: this is to be defined}}%
and thus $\sum_{j = 1}^m {\cal T}_j$ contains any simple spinor $\psi$ and thus $M(\psi)$ can be any MTNP $(\Identity, t)$ proving (\ref{formula_SAT_in_O(n)}). Conversely assuming (\ref{formula_SAT_in_O(n)}) this means that with $\sum_{j = 1}^m {\cal T}_j$ we can express any simple spinor $\psi \in \mySpinorS_s$ and thus any $\psi_\lambda \in \myFockB$ and by Lemma~\ref{lemma_T_T'_intersection}, being $\myFockB$ a basis of $\mySpinorS$, this implies that $\cup_{j = 1}^m {\cal T}_j' = \O1n = \myFockB$ and $\myBooleanS$ is unsatisfiable by Proposition~\ref{SAT_in_n_O(1)}.
\end{proof}

We remark that this theorem is not a straightforward generalization of Proposition~\ref{SAT_in_n_O(1)} and that sum operation (\ref{formula_Tau_addition}) is pivotal: replacing $\sum$ with $\cup$ in (\ref{formula_SAT_in_O(n)}) the result does not hold and there are unsatisfiable problems for which $\cup_{j = 1}^m {\cal T}_j \ne \OO{n}$. For example in previous simple case in $\R^{2,2}$ the \SAT{} problem given by the $4$ clauses $\literal_1 \literal_2, \literal_1 \myconjugate{\literal}_2, \myconjugate{\literal}_1 \literal_2, \myconjugate{\literal}_1 \myconjugate{\literal}_2$ is clearly unsatisfiable but for any $t_\theta \in \SO{2}$ with $\theta \ne 0, \pi$ the isometry $t_\theta \notin \cup_{j = 1}^m {\cal T}_j = \cup_{j = 1}^m {\cal T}_j'$ while, as shown above, $t_\theta \in {\cal T}_{\literal_1 \literal_2} + {\cal T}_{\myconjugate{\literal}_1 \myconjugate{\literal}_2}$.
\opt{margin_notes}{\mynote{mbh.note: log pp. 774, 784}}%

This result gives an unsatisfiability test that examines the clauses $z_j$ to verify if they induce a cover of $\OO{n}$. We just remind that $\OO{n}$ is a continuous group that form a compact, disconnected real manifold of dimension $n (n-1)/2$ and that its two connected components are respectively $\SO{n}$, with $\Identity$, and its coset given by $\OO{n}$ elements with determinant $-1$ \cite{Porteous_1995}.

\section{Conclusions and outlook}
\label{Conclusions}
In the first part of the paper we have shown that \SAT{} fits neatly in Clifford algebra allowing to look at \SAT{} from a different viewpoint. This culminated in Proposition~\ref{SAT_in_n_O(1)} that substantially says that to test unsatisfiability we have to verify if the group \O1n is covered by the involutions induced by its clauses. This implies that any $\lambda \in \O1n$ not in $\cup_{j = 1}^m {\cal T}_j'$ is a solution of the problem at hand.

It is relevant to ask whether the new formulation contributes also to the computational side of \SAT{} and the answer is \emph{no}. We crudely resume the current situation saying that any \SAT{} algorithm looks after a solution, namely an assignment of the $n$ \myBooleanv{}s, and that, in worst case, it will need to check ${\cal O}(2^n)$ assignments to find a solution or to prove that there are none. An algorithm willing to exploit Proposition~\ref{SAT_in_n_O(1)} would search after an element of \O1n not contained in $\cup_{j = 1}^m {\cal T}_j'$. Being \O1n a discrete group essentially the only possibility is to check one by one each of its $2^n$ elements. On top of that it is easy to prove also that testing if any given ${\cal T}_x'$ is contained in $\cup_{j = 1}^m {\cal T}_j'$ is NP-complete%
\footnote{given clause $z_x$ that makes the problem certainly unsatisfiable then ${\cal T}_x' \subseteq \cup_{j = 1}^m {\cal T}_j'$ would provide a certificate of unsatisfiability.}%
. Put in this way the setting in Clifford algebra does not bring any substantial advantage with respect to usual algorithms.

The scene changes in Section~\ref{continuous_isometries}: thanks to pivotal Proposition~\ref{prop_5_BudinichP_1989} we can exploit a unique property of $\R^{n,n}$ applying it to \SAT{} formulation in Clifford algebra. Theorem~\ref{SAT_in_O(n)} shows that an unsatisfiable problem induces also a full cover of $\OO{n}$, a continuous group. Again an algorithm could test unsatisfiability checking if there exist $t \in \OO{n}$ not contained in $\sum_{j = 1}^m {\cal T}_j$ but now the continuity of $\OO{n}$ makes the situation radically different with respect to the case of \O1n.

We give an argument to support this claim: exploiting a parametrization of $\OO{n}$ elements we can transform the sets ${\cal T}$ (\ref{formula_cal_T_j_def2}) induced by clauses into subsets of the parameter space (for example to subsets of $[0, 2 \pi)^{\frac{n(n - 1)}{2}}$ in the case of decomposition in Givens rotations, bivectors in Clifford algebra) to search $\OO{n}$ elements not contained in $\cup_{j = 1}^m {\cal T}_j$ in parameter space. This is more similar to searching if a given simple spinor can be built out of the various subsets of Fock basis $\myFockB$ induced by clauses ${\cal T}$ (\ref{formula_cal_T_j_def2}) and there seems to be no combinatorial calculus in sight. This hints a path to follow that, even if challenging, appears to be a non beaten track heading to unexplored territories and a worthy subject for future research. The same setting seems also suitable to reformulate other \SAT{} theoretical issues like \eg threshold phenomena for random \SAT{} instances.

%
%
%
%
%
%
%Being $\OO{n}$ a smooth, low dimensional, manifold its exploration could in principle take advantage from the deployment of the full arsenal of tools of analysis.
%
%The difference is that while \O1n is a discrete group and to verify (\ref{formula_SAT_in_n_O(1)}) one essentially needs to check one by one all possible $2^n$ transformations, $\OO{n}$ is a compact and disconnected real manifold of dimension $n (n-1)/2$ \cite{Porteous_1995} that could make simpler to verify -- thus checking $S$ for unsatisfiability -- if the union of the isometries induced by the set of clauses $z_j$ form a cover of $\OO{n}$.

%\bigskip
%\bigskip
%\bigskip
%\bigskip
%\bigskip

\newpage

\opt{x,std,AACA}{% in all cases - but arXiv, JMP & TCS - standard BibTeX bibliography

\bibliographystyle{plain} % or plain or.... see e.g.\ http://amath.colorado.edu/documentation/LaTeX/reference/faq/bibstyles.html#styles
\bibliography{mbh}
}

\opt{arXiv,JMP,TCS}{% only for arXiv, JMP & TCS we need to include here the L-A-S-T version of file .bbl
%
%
%\begin{thebibliography}{1}

%\end{thebibliography}
%
%
}

\opt{final_notes}{

\newpage

\section*{Appendix: on Witt bases of $\R^{n,n}$ (OLD, paper I vrs)}
\label{Appendix}
A Witt basis of $\R^{n,n}$ is defined by two MTNP $P$ and $Q$ of \mysnqG{} such that $P \oplus Q = \R^{n,n}$ and since by (\ref{formula_bijection_Nn_O(n)}) any MTNP can be represented in the form $(\Identity, t)$ with $t \in \OO{n}$ we will indicate the Witt basis as column vectors by $W = (P, Q) = \left(\begin{array}{r r} \Identity & \Identity \\ t_1 & t_2 \end{array}\right)$.

\begin{MS_lemma}
\label{MTNP_equivalence}
Given any $r, s \in GL(n)$ such that $s r^{-1} \in \OO{n}$ then $(r, s)$ and $(\Identity, s r^{-1})$ represent the same MTNP.
\opt{margin_notes}{\mynote{mbh.note: very similar to \cite[Proposition~12.14]{Porteous_1981}, for a fully general version see log. p. 721'}}%
\end{MS_lemma}
\begin{proof}
Putting the two subspaces in a $2n \times 2n$ matrix it is easy to verify that
$$
\left(\begin{array}{c c} r & \Identity \\ s & s r^{-1} \end{array}\right) \left(\begin{array}{c} a \\ - r (a) \end{array}\right) = \left(\begin{array}{c} 0 \\ 0 \end{array}\right) \qquad \forall a \in \R^n
$$
and thus the matrix has rank $n$ since it has a null space of dimension $n$ and $(\Identity, s r^{-1})$ is a MTNP and thus has rank $n$.
\end{proof}
\noindent From now on we will take any MTNP to be of the form $(\Identity, t)$ with $t \in \OO{n}$.

\begin{MS_Proposition}
\label{condition_Witt_basis}
Given two MTNP $(\Identity, t_1)$ and $(\Identity, t_2)$ they form a Witt basis if and only if $\det(t_2 - t_1) \ne 0$.
\end{MS_Proposition}
\begin{proof}
The two MTNP form a Witt basis if and only if $\det \left(\begin{array}{c c} \Identity & \Identity \\ t_1 & t_2 \end{array}\right) \ne 0$ but by block matrices properties this determinant is $\det(t_2 - t_1)$.
%
%\opt{margin_notes}{\mynote{mbh.note: $\det \left(\begin{array}{r r} A & B \\C & D \end{array}\right) = \det (A) \det (D - C A^{-1} B)$}}%
%
\end{proof}
\noindent It follows that MTNP $P = (\Identity, \Identity)$ and $(\Identity, t)$ form a Witt basis if and only if $\det(t - \Identity) \ne 0$; since the eigenvalues of all elements of $\OO{n}$ belong to $\{ \pm 1, e^{\pm i \theta} \}$ we deduce that for our $t \in \OO{n}$ there are not eigenvalues $1$.

For Witt bases we can show a property similar to the unicity of the orthogonal complement of a non degenerate subspace in Euclidean spaces:
\begin{MS_Proposition}
\label{complementary_MTNP_unicity}
Given any MTNP $P \in \mysnqG{}$ there always exists one and only one complementary MTNP $Q_\perp$ such that $P \oplus Q_\perp = \R^{n,n}$.
\opt{margin_notes}{\mynote{mbh.note: does this show $Q$ in a basis independent form ?}}%
\end{MS_Proposition}

\begin{proof}
Let $P := (\Identity, t)$ with $t \in \OO{n}$; since $\frac{1}{2} P^T P = \Identity_n$ $P$ represents also a non degenerate subspace of dimension $n$ of Euclidean space $R^{2 n}$. Let $Q := (A, B)$ be the most general linear complement of $P$ in $R^{2 n}$, then
$$
(P, Q)^T (P, Q) = \left(\begin{array}{r r} \Identity & t^T \\ A^T & B^T \end{array}\right) \left(\begin{array}{r r} \Identity & A \\ t & B \end{array}\right) = \left(\begin{array}{c c} 2 \,\Identity & A + t^T B \\ A^T + B^T t & A^T A + B^T B \end{array}\right) \dotinformula
$$
Among the linear complements of $P$ there exists a unique orthogonal complement $Q_\perp$ \cite[Proposition~9.25]{Porteous_1981} that must have $B = -t A$ and has thus the form $Q_\perp := (A, -t A)$ so that in $R^{2 n}$, with proper normalization,
$$
\frac{1}{2} (P, Q_\perp)^T (P, Q_\perp) = \frac{1}{2} \left(\begin{array}{c c} \Identity & t^T \\ A^T & - A^T t^T \end{array}\right) \left(\begin{array}{c c} \Identity & A \\ t & -t A \end{array}\right) = \left(\begin{array}{c c} \Identity & 0 \\ 0 & A^T A \end{array}\right)
$$
\opt{margin_notes}{\mynote{mbh.note: in $\R^{n,n}$ this gives a non orthogonal Witt basis with $A$ \& $A^T$ at place of $\Identity$%
%$\left(\begin{array}{r r} 0 & A \\ A^T & 0 \end{array}\right)$%
}}%
that shows that $A^T A$, and thus $A$, must be full rank and so $A \in GL(n)$.
%
%Getting back to $\R^{n,n}$ for the same subspaces we get
%$$
%\frac{1}{2} (P, Q_\perp)^T \left(\begin{array}{r r} \Identity & 0 \\ 0 & - \Identity \end{array}\right) (P, Q_\perp) = \frac{1}{2} \left(\begin{array}{c c} \Identity & -t^T \\ A^T & A^T t^T \end{array}\right) \left(\begin{array}{c c} \Identity & A \\ t & -t A \end{array}\right) = \left(\begin{array}{c c} 0 & A \\ A^T & 0 \end{array}\right)
%$$
%and thus they form a Witt basis since $A$ is full rank. In $R^{2 n}$ the unique orthogonal complement of $P$ has the form $Q_\perp = (A, -t A)$ with $A$ full rank and thus
%
We can thus rearrange the basis in $Q_\perp$ to obtain $A = \Identity_n$ so in this basis $Q_\perp = (\Identity, -t)$. Back in $\R^{n,n}$ it is simple to verify that MTNP $P$ and $Q_\perp = (\Identity, -t)$ form an orthogonal Witt basis (\ref{formula_Witt_orthogonality_def})
\begin{equation}
\label{formula_Witt_orthogonality}
\frac{1}{2} \left(\begin{array}{r r} \Identity & \Identity \\ t & - t \end{array}\right)^T \left(\begin{array}{r r} \Identity & 0 \\ 0 & - \Identity \end{array}\right) \left(\begin{array}{r r} \Identity & \Identity \\ t & - t \end{array}\right) = \left(\begin{array}{r r} 0 & \Identity \\ \Identity & 0 \end{array}\right) \dotinformula
\end{equation}
Should there exist another, different, MTNP $(\Identity, t')$ forming a Witt basis with $P$, $(\Identity, t')$ would violate the unicity of the orthogonal complement of $P$ in $R^{2 n}$ thus the unicity of MTNP $Q_\perp$.
\opt{margin_notes}{\mynote{mbh.note: is the proof deducible immediately from Witt cancellation theorem ?}}%
\end{proof}
In summary given MTNP subspace $P$ of $\R^{n,n}$ it has very many linear complements but the additional request that the linear complement is also a MTNP makes the solution unique, this being identical to what happens in Euclidean space: given a non degenerate subspace there are many linear complements but only one of them is also an orthogonal complement.

We define a Witt basis $W$ of $\R^{n,n}$ to be orthogonal iff
\begin{equation}
\label{formula_Witt_orthogonality_def}
W^T \left(\begin{array}{r r} \Identity & 0 \\ 0 & - \Identity \end{array}\right) W = \left(\begin{array}{r r} 0 & \Identity \\ \Identity & 0 \end{array}\right)
\end{equation}
\opt{margin_notes}{\mynote{mbh.note: this form is not immediate for the basis definition (\ref{formula_generators}), some comments should be added}}%
for example the standard Witt basis $\frac{1}{\sqrt{2}} \left(\begin{array}{r r} \Identity & \Identity \\ \Identity & - \Identity \end{array}\right) := W_0$ is orthogonal; we remark that $W_0$ is essentially the Witt basis of (\ref{formula_Witt_basis}) with a normalization slightly different from that usually adopted in EFB that turns out to be more handy when dealing with linear spaces.
\opt{margin_notes}{\mynote{mbh.note: the normalization of (\ref{formula_Witt_basis}) and (\ref{formula_P_Q_def}) is treated in logbook pp. 343 ff.}}%
With this definition:
\begin{MS_Corollary}
\label{complementary_MTNP_unicity_coro}
Given any MTNP $P \in \mysnqG{}$ with $P = (\Identity, t)$ its unique MTNP complement $Q$ can take all possible forms given by $(A, - t A)$ for any $A \in GL(n)$ and the Witt basis $(P, Q)$ results orthogonal if and only if $A = \Identity_n$.
\end{MS_Corollary}

\begin{MS_Proposition}
\label{orthogonal_Witt_basis}
A Witt basis $W$ is orthogonal if and only if it is of the form $W = \frac{1}{\sqrt{2}} \left(\begin{array}{r r} \Identity & \Identity \\ t & -t \end{array}\right)$ for any $t \in \OO{n}$.
\end{MS_Proposition}
\begin{proof}
We already saw (\ref{formula_Witt_orthogonality}) that $W$ is orthogonal, conversely let $W = \frac{1}{\sqrt{2}} \left(\begin{array}{r r} \Identity & \Identity \\ t_1 & t_2 \end{array}\right)$ be an orthogonal Witt basis then
$$
\left(\begin{array}{r r} 0 & \Identity \\ \Identity & 0 \end{array}\right) = \frac{1}{2} \left(\begin{array}{r r} \Identity & \Identity \\ t_1 & t_2 \end{array}\right)^T \left(\begin{array}{r r} \Identity & 0 \\ 0 & - \Identity \end{array}\right) \left(\begin{array}{r r} \Identity & \Identity \\ t_1 & t_2 \end{array}\right) = \frac{1}{2} \left(\begin{array}{c c} 0 & \Identity - t_1^T t_2 \\ \Identity - t_2^T t_1 & 0 \end{array}\right)
$$
and thus $\frac{1}{2} (\Identity - t_1^T t_2) = \frac{1}{2} (\Identity - t_2^T t_1) = \Identity$, namely $t_2 = - t_1$.
\end{proof}
\noindent This shows that bijection (\ref{formula_bijection_Nn_O(n)}) extends also to orthogonal Witt bases $W$; moreover it is immediate to verify:
\opt{margin_notes}{\mynote{mbh.note: how do we express a cover of $\OO{n}$ in terms of Witt bases ?}}%
\begin{MS_Corollary}
\label{Witt_is_on_basis}
Any orthogonal Witt basis of $\R^{n,n}$ is also an orthonormal basis of $\R^{2 n}$ since $W^T W = \Identity$; conversely any orthonormal basis of $\R^{2 n}$ of the form $W = \frac{1}{\sqrt{2}} \left(\begin{array}{c c} t_1 & t_1 \\ t_2 & - t_2 \end{array}\right)$ with $t_1, t_2 \in \OO{n}$, is also an orthogonal Witt basis of $\R^{n,n}$.
\end{MS_Corollary}

We conclude with an alternative proof of (a restricted form of) the pivotal Proposition~\ref{prop_2_BudinichP_1989} without spinors but with standard matrix formalism.

\begin{MS_Proposition}
\label{prop_2_BudinichP_1989_matricial}
Given any Witt basis $W = \frac{1}{\sqrt{2}} \left(\begin{array}{r r} \Identity & \Identity \\ \Identity & t \end{array}\right)$ of $\R^{n,n}$ then there always exists an orthonormal basis of $\R^{n,n}$ in which the two MTNP forming $W$ assume the form $W_0$, in particular MTNP $(\Identity, \Identity)$ keeps its form while $(\Identity, t)$ takes the form $(\Identity, - \Identity)$.
\end{MS_Proposition}

\begin{proof}
Since $W$ is a Witt basis the matrix $W$ is invertible and so (non orthogonal) rotation $W_0 W^{-1}$ transforms the null vectors forming $W$ to the orthogonal Witt basis $W_0$. We can calculate it explicitly: since $W$ is a Witt basis by Proposition~\ref{condition_Witt_basis} $\det(t - \Identity) \ne 0$ and thus $(t - \Identity)^{-1}$ exists; it is easy to verify that
$$
W^{-1} = \left(\begin{array}{c c} (t - \Identity)^{-1} t & - (t - \Identity)^{-1} \\ - (t - \Identity)^{-1} & (t - \Identity)^{-1} \end{array}\right)
$$
\opt{margin_notes}{\mynote{mbh.note: remember $(\Identity - t)^{-1} = \sum_{k = 0}^\infty t^k$}}%
and thus the required rotation is
$$
W_0 W^{-1} = \frac{1}{\sqrt{2}} \left(\begin{array}{c c} \Identity & 0 \\ (t - \Identity)^{-1} (t + \Identity) & -2 (t - \Identity)^{-1} \end{array}\right)
$$
that leaves invariant $(\Identity, \Identity)$ all its null vectors being eigenvectors of $W_0 W^{-1}$ of eigenvalue $\frac{1}{\sqrt{2}}$.
\opt{margin_notes}{\mynote{mbh.ref: the MTNP complement of $(\Identity, \Identity)$ is unique: see log. p. 717.}}%
By Proposition~\ref{complementary_MTNP_unicity} the unique MTNP complement of $(\Identity, \Identity)$ in $\R^{n,n}$ has been thus transformed by $W_0 W^{-1}$ from $(\Identity, t)$ to $(\Identity, - \Identity)$.
\end{proof}
\noindent We remark that all $t \in \OO{n}$ with all eigenvalues $\ne 1$ represent MTNP $(\Identity, - \Identity)$.

%\bigskip

\begin{MS_lemma}
\label{MTNP_reference}
Given any two MTNP forming $W = \left(\begin{array}{c c} \Identity & \Identity \\ t_1 & t_2 \end{array}\right)$ it is always possible to rearrange the basis in $\{0\} \times \R^n$ so that the null vectors forming $W$ in the new basis take the form $W' = \left(\begin{array}{c c} \Identity & \Identity \\ \Identity & t \end{array}\right)$ with $t = t_1^T t_2$.
\end{MS_lemma}
\begin{proof}
Given orthonormal basis $E = \left(\begin{array}{c c} \Identity & 0 \\ 0 & t_1 \end{array}\right)$ then $W' = E^{-1} W = \left(\begin{array}{c c} \Identity & \Identity \\ \Identity & t_1^T t_2 \end{array}\right)$.
\end{proof}

\newpage

\section*{Things to do, notes, etc.......}

%\newpage
\section{\SAT{} Clauses and continuous isometries of $\OO{n}$ (Old version of paragraph 6 of VI 2021)}
\label{continuous_isometries_old}
In the last step we show that when a \SAT{} problem $S$ is unsatisfiable the isometries induced by its clauses not only form a cover of \O1n (\ref{formula_SAT_in_n_O(1)}) but also of its parent group $\OO{n}$.

To achieve this we need to update the definition of isometries induced by a clause (\ref{formula_cal_T_j_def}) in a basis independent form and taking advantage of a special property of Witt bases of $\R^{n,n}$.

Proposition~\ref{isomorphism_restricted} establishes a one to one correspondence between $\lambda \in \O1n$ and \mysetM{} but changing basis of $\R^{n,n}$ the same MTNP is associated to a possibly different $\lambda'$ and by bijection (\ref{formula_bijection_Nn_O(n)}) $\lambda' \in \OO{n}$. We will thus associate to a MTNP of \mysetM{} not only $\lambda \in \O1n$ but the whole set of $t \in \OO{n}$ representing it in all possible forms passing from ${\cal T}_j'$ of (\ref{formula_cal_T_j_def}) to a superset ${\cal T}_j \supset {\cal T}_j'$ that widens the definition of the set of isometries induced by $z_j$
\opt{margin_notes}{\mynote{mbh.note: tricky point! also $M(z_j)$ could change its form. Probably it shouldn't be a problem since there are only $2^n$ MTNP; moreover the ``form'' of $(\Identity, t)$ doesn't depend only on basis}}%
to
\begin{equation}
\label{formula_cal_T_j_def2_old}
{\cal T}_j := \{t_j \in \OO{n} : M(z_j) \subseteq (\Identity, t_j)% \; \mbox{in all bases of} \; \R^{n,n}
\} \supset {\cal T}_j' \dotinformula
\end{equation}
We remark that since ${\cal T}_j'$ contains in general $2^{n-k}$ isometries of $\O1n$ the corresponding superset ${\cal T}_j$ will contain all $t \in \OO{n}$ representing the $2^{n-k}$ isometries of $\O1n$ in all possible forms.

To proceed we need a property of Witt bases \cite{BudinichP_1989} that we reproduce here, slightly adapted to our needs:
\begin{MS_Proposition}
\label{prop_2_BudinichP_1989}
Given any two MTNP of \mysnqG{} $(\Identity, t_1)$ and $(\Identity, t_2)$ with $t_1 \ne t_2$ then it is always possible to choose an orthonormal basis of $\R^{n,n}$ such that the two given MTNP in the new basis take the form of two MTNP of \mysetM{}.
\end{MS_Proposition}
Any MTNP of \mysnqG{} is a subspace of $\R^{n,n}$ and two of them necessarily have an intersection of dimension $r$ with $0 \le r \le n$. For all cases of $r < n$ one can choose a basis of $\R^{n,n}$ such that the two given MTNP coincide with two MTNP of \mysetM{} \cite{BudinichP_1989}. This can be rephrased saying that given two different MTNP of \mysnqG{} one can always find a basis of $\R^{n,n}$ in which the two MTNP take the form of $P = (\Identity, \Identity)$ and $(\Identity, \lambda)$ with $\lambda \in \O1n$.

For example in $\R^{2,2}$ given MTNP $P = (\Identity, \Identity)$ and $(\Identity, t_\theta)$, where $t_\theta \in \SO{2}$ with angle $\theta \ne 0$, we can choose a basis so that the very same MTNP $(\Identity, t_\theta)$ is represented by $Q = (\Identity, - \Identity)$ in the new basis and what matters is not the form $(\Identity, t_\theta)$ of the MTNP in a particular basis but only the fact that its intersection with $P$ is $\{ 0 \}$ and thus that their direct sum is $\R^{n,n}$.
\opt{margin_notes}{\mynote{mbh.note: probably one should prove this statement in $\R^{2,2}$ for example showing that $P \cap (\Identity, t_\theta) = \{0\}$ for all $\theta \ne 0$}}%
We remark that this property holds only for null bases of $\R^{n,n}$ that have the peculiar property that any invertible linear combination (not necessarily orthogonal) of the basis of $P$ produces a new basis of vectors that remain reciprocally orthogonal, a property not holding \eg in Euclidean spaces.

\begin{MS_theorem}
\label{SAT_in_O(n)_old}
A given \SAT{} problem $S$ with m clauses $z_j$ is unsatisfiable if and only if the isometries induced by its clauses (\ref{formula_cal_T_j_def2_old}) form a cover for $\OO{n}$:
\begin{equation}
\label{formula_SAT_in_O(n)_old}
\cup_{j = 1}^m {\cal T}_j = \OO{n} \dotinformula
\end{equation}
\end{MS_theorem}
\begin{proof}
Let (\ref{formula_SAT_in_O(n)_old}) hold, since ${\cal T}_j' \subset {\cal T}_j$ if $\cup_{j = 1}^m {\cal T}_j' = \O1n$ then $S$ is unsatisfiable by Proposition~\ref{SAT_in_n_O(1)}. To achieve this we need to prove that for any $j$ there are no $\lambda \in \O1n$ in ${\cal T}_j$ that are not in ${\cal T}_j'$, namely that ${\cal T}_j \cap \O1n = {\cal T}_j'$. We prove this observing that for any $\lambda \in \O1n - {\cal T}_j'$ then $M(z_j) \nsubseteq (\Identity, \lambda)$ (\ref{formula_cal_T_j_def}) and thus, in all its possible forms, including $\lambda$ itself, it cannot be in ${\cal T}_j$ (\ref{formula_cal_T_j_def2_old}) and so $\cup_{j = 1}^m {\cal T}_j' = \O1n$ and $S$ is unsatisfiable.

Conversely let $S$ be unsatisfiable, for $n = 1$ $\O11 = \mbox{O}(1)$ and (\ref{formula_SAT_in_O(n)_old}) and (\ref{formula_SAT_in_n_O(1)}) coincide so there is nothing to prove.
\opt{margin_notes}{\mynote{mbh.note: by Proposition~\ref{M_z_j} we need $k < n-2$ so perhaps one needs an explicit proof for the case $n = 2$; commented here there is this proof}}%
For $n > 1$ by Proposition~\ref{SAT_in_n_O(1)} and (\ref{formula_cal_T_j_def2_old}) the union (\ref{formula_SAT_in_O(n)_old}) covers $\O1n$ and thus contains the $2^n$ involutions $\lambda_j \in \O1n$ corresponding to the MTNP of \mysetM{} and also all the isometries $t_j \in \OO{n}$ representing these MTNP in any possible form. It remains to be proved that any $t \in \OO{n}$ is in this set. Let us suppose the contrary, namely that there exists $t' \in \OO{n}$ that is not in this set: by bijection (\ref{formula_bijection_Nn_O(n)}) $(\Identity, t')$ is a MTNP and a subspace of $\R^{n,n}$ that necessarily must have an intersection of dimension $0 \le r \le n$ with our reference MTNP $P$ and then by Proposition~\ref{prop_2_BudinichP_1989} there exists a basis of $\R^{n,n}$ such that this MTNP is in \mysetM{} and thus in the union (\ref{formula_SAT_in_O(n)_old}) contradicting our hypothesis. We must necessarily conclude that any $t \in \OO{n}$ is in $\cup_{j = 1}^m {\cal T}_j$.
\end{proof}

This result gives an unsatisfiability test that examines the clauses $z_j$ to verify if they induce a cover of $\OO{n}$. We just remind that $\OO{n}$ is a continuous group that form a compact, disconnected real manifold of dimension $n (n-1)/2$ and that the two connected components are respectively $\SO{n}$, with $\Identity$, and its coset given by $\OO{n}$ elements with determinant $-1$ \cite{Porteous_1995}.

\subsection*{Other Propositions, other material, etc.......}

We already remarked that there is a one to one correspondence between \mysnqG{} and $\OO{n}$ by (\ref{formula_bijection_Nn_O(n)}) and thus since for our reference MTNP $P = (\Identity, \Identity)$ we can see any MTNP as the 'action' of $t(x)$ on $P$ by means of $(\Identity, t \circ \Identity(x)) = (\Identity, t(x))$. What we are going to show is that while the explicit form of $t \in \OO{n}$ depends on the basis of $\R^{n,n}$ we can define any MTNP, or its $t \in \OO{n}$, in a basis independent form and this allows to associate to each MTNP a set of equivalent $t \in \OO{n}$.
\opt{margin_notes}{\mynote{mbh.note: this shows that the two possible defs of $t$, as a MTNP or as an action on a MTNP are equivalent}}%

\bigskip

Possible other roads to explore, good questions etc:
\begin{itemize}
\item induction on $n$ (see below a tentative proof in this direction).
\item We can build all the $n$ subgroups $\mbox{O}(1) \times \OO{n-1}$.
\item Given $\OO{n}$ it has $2^n$ subgroups $\OO{k} \times \OO{n-k}$; is $\cup_{i = 1}^{2^n} \OO{k} \times \OO{n-k} = \OO{n}$ ? If this is true we have just to prove that if a problem is unsat then each ${\cal T}_j$ contains $2^k + 2^{n-k}$ subgroups and that their union covers $\cup_{i = 1}^{2^n} \OO{k} \times \OO{n-k}$ that probably is simpler to prove.
\item The most general element of \mysnqG{} can be written either in the form $(x, t (x))$ or $(u(x), t u (x))$.
\item It is NOT true that all elements of \mysnqG{} of the form $(x, t(x))$ correspond to one of the $2^n$ MTNP. Counterexample in $\R^{2,2}$:
$$
\my_span{e_1 + \frac{1}{\sqrt{2}} (e_3 - e_4), e_2 + \frac{1}{\sqrt{2}} (e_3 + e_4)}
$$
is a MTNP given by $t_\theta = \frac{1}{\sqrt{2}} \left(\begin{array}{r r} 1 & -1 \\ 1 & 1 \end{array}\right) = t_{\frac{\pi}{4}} \in \OO{2}$. What IS true is that exists a basis of $\R^{n,n}$ in which they take the form $\left(\begin{array}{r r} \Identity & \Identity \\ \Identity & \lambda \end{array}\right)$ with $\lambda \in \O1n$ and that $\det() \ne 0$ iff $\lambda = - \Identity$.
\end{itemize}

\begin{proof}
\opt{margin_notes}{\mynote{mbh.note: this is an alternative, correct, proof of Theorem~\ref{SAT_in_O(n)_old}}}%
OLD PROOF n. 4 - Let $S$ be unsatisfiable, by Proposition~\ref{SAT_in_n_O(1)} we know that (\ref{formula_SAT_in_n_O(1)}) holds and thus that $\O1n \subset \cup_{j = 1}^m {\cal T}_j$. By Proposition~\ref{MTNPs} any $t \in \OO{n}$ corresponds to one element of \mysetM{} in a particular basis of $\R^{n,n}$ and thus is in at least one ${\cal T}_j$ that proves (\ref{formula_SAT_in_O(n)_old}).

%Exploiting the isomorphism between between $\OO{n}$ and \mysnqG{} to prove (\ref{formula_SAT_in_O(n)_old}) it is sufficient to show that $\cup_{j = 1}^m {\cal T}_j$ acts transitively on \mysnqG{}. By Proposition~\ref{MTNPs} any $t \in \OO{n}$ gives one element of \mysetM{} that are already in $\cup_{j = 1}^m {\cal T}_j$.
%
%Any isometry $t \in \OO{n}$ identify a MTNP of \mysnqG{} that, by Proposition~\ref{MTNPs} is one of the possible $2^n$ MTNP that proves that $t \in \cup_{j = 1}^m {\cal T}_j$ and thus (\ref{formula_SAT_in_O(n)_old}).

Conversely ....
\end{proof}

\begin{proof}
OLD PROOF n. 2 - Let $S$ be unsatisfiable, we start proving that any MTNP of \mysnqG{} identified by $t \in \OO{n}$ can be written in the form $T = T_k \Lambda_j$ with $\Lambda_j \in \O1n$ and $T_k \in \OO{k} \times \OO{n-k}$ for some $k$ and thus that any $t \in \OO{n}$ may be written as $t = \Lambda_j T_k$ and thus satisfies (\ref{formula_Lambda_j_action}) and is thus one of the isometries induced by $z_j$ (here something is missing).
Let us suppose that $(x, t(x)) \in \mysnqG$ is contained in a MTNP formed by $k$ vectors $q_i$ and $n-k$ $p_i$, then it is possible to write $T = \Lambda_j T_k$ where $\Lambda_j$ selects one of the $2^n$ MTNP of the Witt basis and $T_k \in \OO{k} \times \OO{n-k}$ perform two separate isometries in the span of the $q_i$ and $p_i$. Since \mysnqG{} is isomorphic to $\OO{n}$ this means that we can write any $T \in \OO{n}$ in the form $T = \Lambda_j T_k$ and thus we have shown that the unsatisfiability of $S$ implies (\ref{formula_SAT_in_O(n)_old}).

Conversely ....
\end{proof}

Some previous propositions, correct but not used anymore:

\begin{MS_Proposition}
\label{prop_2_BudinichP_1989_matricial_2}
Given any Witt basis $W = \frac{1}{\sqrt{2}} \left(\begin{array}{r r} \Identity & \Identity \\ \Identity & t \end{array}\right)$ of $\R^{n,n}$ then there exists an orthonormal basis of $\R^{n,n}$ in which the form of the Witt basis $W$ becomes orthogonal with
$$
W' = \left(\begin{array}{c c} \Identity & \Identity \\ \lambda (\Identity - t)^{-1} & - \lambda (\Identity - t)^{-1} \end{array}\right)
$$
(wrong!!) with $\lambda (\Identity - t)^{-1} \in \OO{n}$ and $\lambda \in \R^* := \R - \{0\}$.
\end{MS_Proposition}

\begin{proof}
The matrix $g_n$ of the scalar products in $\R^{n,n}$ of the column vectors forming $W$ is
$$
g_n = \frac{1}{2} W^T \left(\begin{array}{r r} \Identity & 0 \\ 0 & - \Identity \end{array}\right) W = \frac{1}{2} \left(\begin{array}{c c} 0 & \Identity - t \\ \Identity - t^T & 0 \end{array}\right) := \left(\begin{array}{c c} 0 & \alpha \\ \alpha^T & 0 \end{array}\right)
$$
and since by hypothesis $W$ is a Witt basis by Proposition~\ref{condition_Witt_basis} $\det(\alpha) = (-\frac{1}{2})^n \det(t - \Identity) \ne 0$ that proves that $\alpha^{-1}$ exists. $g_n$ is symmetric and thus is diagonalizable and with eigenvalues $\lambda_i \in \R$, moreover since
$$
\det(g_n) = \det( - \alpha^T \alpha) = (-1)^n \det(\alpha)^2 \ne 0
$$
it follows that the eigenvalues $\lambda_i \in \R^*$. There is thus an orthogonal matrix $R$ of the Euclidean space $\R^{2 n}$ and a diagonal matrix $\Lambda_g$ such that
$$
R^T g_n R = \Lambda_g \qquad \implies \qquad g_n R = R \Lambda_g \dotinformula
$$
Looking at $\R^{2 n}$ as $\R^{n} \times \R^{n}$ any eigenvector of $R$ can be written in the form $(a,b)$ and thus we must have
$$
g_n \left(\begin{array}{c} a \\ b \end{array}\right) = \left(\begin{array}{c c} 0 & \alpha \\ \alpha^T & 0 \end{array}\right) \left(\begin{array}{c} a \\ b \end{array}\right) = \left(\begin{array}{c} \alpha b \\ \alpha^T a \end{array}\right) = \lambda \left(\begin{array}{c} a \\ b \end{array}\right)
$$
that is satisfied only if
\begin{equation}
\label{formula_a_b_conditions}
\left\{ \begin{array}{l l l}
\alpha b & = & \lambda a \\
\alpha^T a & = & \lambda b \dotinformula
\end{array} \right.
\end{equation}
Before going on we remark that if $(a,b)$ is an eigenvector of $g_n$ of eigenvalue $\lambda$, and thus (\ref{formula_a_b_conditions}) are satisfied, then
$$
g_n \left(\begin{array}{c} a \\ - b \end{array}\right)
%= \left(\begin{array}{c c} 0 & \alpha \\ \alpha^T & 0 \end{array}\right) \left(\begin{array}{c} a \\ - b \end{array}\right)
= \left(\begin{array}{c} - \alpha b \\ \alpha^T a \end{array}\right) = \left(\begin{array}{c} - \lambda a \\ \lambda b \end{array}\right) = - \lambda \left(\begin{array}{c} a \\ - b \end{array}\right)
$$
and thus $(a,-b)$ is an eigenvector of $g_n$ of eigenvalue $- \lambda$ that proves that all $2 n$ eigenvalues of $g_n$ are in couples $\pm \lambda_i$ and thus we can write $\Lambda_g = \left(\begin{array}{r r} \Lambda & \\ & - \Lambda \end{array}\right)$ where we can rearrange the eigenvalues in such a way that all $n$ eigenvalues in $\Lambda$ are positive. Getting back to the solution of (\ref{formula_a_b_conditions}) since $\lambda \ne 0$ we get that any $a \in \R^n$, part of an eigenvector of $g_n$, must satisfy
\opt{margin_notes}{\mynote{mbh.note: alternatively $\alpha^T \alpha b = \lambda^2 b$}}%
\begin{equation}
\label{formula_a_b_conditions_2}
\alpha \alpha^T a = \lambda^2 a
\end{equation}
and thus that $a$ is an eigenvector of $\alpha \alpha^T$ of eigenvalue $\lambda^2 > 0$. Since the $n \times n$ matrix $\alpha \alpha^T$ is symmetric it has $n$ orthogonal eigenvectors $a_i \in \R^n$ that we put in an orthogonal matrix $A$ and thus
$$
A^T \alpha \alpha^T A = \Lambda^2 \qquad \implies \qquad \alpha \alpha^T A = A \Lambda^2
$$
where $\Lambda^2$ is the square of $\Lambda$ appearing in $\Lambda_g$; from this relation we get immediately that the matrix $\alpha^T A \Lambda^{-1}$, where clearly $\Lambda^{-1} \Lambda = \Identity$, is an orthogonal matrix and thus an element of $\OO{n}$. We can write the eigenvector of $g_n$ associated to eigenvalue $\lambda_i$ and with (\ref{formula_a_b_conditions_2})
\opt{margin_notes}{\mynote{mbh.note: alternatively $\lambda_i \alpha^{-1} a_i$ (below)}}%
$$
\left(\begin{array}{c c} 0 & \alpha \\ \alpha^T & 0 \end{array}\right) \left(\begin{array}{c} a_i \\ \lambda_i^{-1} \alpha^T a_i \end{array}\right) = \left(\begin{array}{c} \lambda_i^{-1} \alpha \alpha^T a_i \\ \alpha^T a_i \end{array}\right) = \left(\begin{array}{c} \lambda_i a_i \\ \alpha^T a_i \end{array}\right) = \lambda_i \left(\begin{array}{c} a_i \\ \lambda_i^{-1} \alpha^T a_i \end{array}\right) \dotinformula
$$
Moreover these relations hold for all $n$ $a_i$ and so we can write the matrix $R$ of the $2 n$ eigenvectors of $g_n$
$$
R = \frac{1}{\sqrt{2}} \left(\begin{array}{c c} A & A \\ \alpha^T A \Lambda^{-1} & - \alpha^T A \Lambda^{-1} \end{array}\right) := \frac{1}{\sqrt{2}} \left(\begin{array}{c c} A & A \\ B & - B \end{array}\right)
$$
with $A$ and $B$ both orthogonal and in $\OO{n}$ and thus
$$
R^T R = \frac{1}{2} \left(\begin{array}{c c} A^T & B^T \\ A^T & - B^T \end{array}\right) \left(\begin{array}{c c} A & A \\ B & - B \end{array}\right) = \left(\begin{array}{c c} \Identity & 0 \\ 0 & \Identity \end{array}\right)
$$
%$$
%R^T R = \left(\begin{array}{c c} A^T & \Lambda^{-1} A^T \alpha \\ A^T & - \Lambda^{-1} A^T \end{array}\right) \left(\begin{array}{c c} A & A \\ \alpha^T A \Lambda^{-1} & - \alpha^T A \Lambda^{-1} \end{array}\right) =
%$$
that confirms that $R$ is the orthogonal rotation of $\OO{2 n}$ that diagonalizes $g_n$; by Corollary~\ref{Witt_is_on_basis} $R$ is also an orthogonal Witt basis of $\R^{n,n}$.
\end{proof}

A different form of previous Proposition.

\begin{MS_Proposition}
\label{prop_2_BudinichP_1989_matricial_3}
Given any Witt basis $W = \frac{1}{\sqrt{2}} \left(\begin{array}{r r} \Identity & \Identity \\ \Identity & t \end{array}\right)$ of $\R^{n,n}$ then ...
\end{MS_Proposition}

\begin{proof}
By hypothesis $W$ is a Witt basis and by Proposition~\ref{condition_Witt_basis} $\det(W) = 2^{-n} \det(t - \Identity) \ne 0$; we easily find the matrix of the scalar products in the Euclidean space $\R^{2 n}$ of the column vectors forming $W$
$$
W^T W = \left(\begin{array}{c c} \Identity & \frac{\Identity + t}{2} \\ \frac{\Identity + t^T}{2} & \Identity \end{array}\right) := \left(\begin{array}{c c} \Identity & \beta \\ \beta^T & \Identity \end{array}\right)
$$
that is symmetric and thus is diagonalizable with eigenvalues $\lambda_i \in \R$, moreover since $\det(W^T W) = \det(W)^2 \ne 0$ and since $W^T W$ is positive definite we can conclude that $\lambda_i > 0$. There is thus an orthogonal matrix $P$ of $\R^{2 n}$ and a diagonal matrix $\Lambda_{W^2}$ such that
$$
P^T W^T W P = \Lambda_{W^2} \qquad \implies \qquad W^T W P = P \Lambda_{W^2} \dotinformula
$$
Looking at $\R^{2 n}$ as $\R^{n} \times \R^{n}$ any eigenvector of $P$ can be written in the form $(a,b)$ and thus we must have
$$
W^T W \left(\begin{array}{c} a \\ b \end{array}\right) = \left(\begin{array}{c c} \Identity & \beta \\ \beta^T & \Identity \end{array}\right) \left(\begin{array}{c} a \\ b \end{array}\right) = \left(\begin{array}{c} a + \beta b \\ \beta^T a + b \end{array}\right) = \lambda \left(\begin{array}{c} a \\ b \end{array}\right)
$$
that is satisfied only if
\begin{equation}
\label{formula_a_b_conditions_bis}
\left\{ \begin{array}{l l l}
\beta b & = & (\lambda - 1) a \\
(\lambda - 1) b & = & \beta^T a \dotinformula
\end{array} \right.
\end{equation}
Before going on we remark that if $(a,b)$ is an eigenvector of $W^T W$ of eigenvalue $\lambda$, and thus (\ref{formula_a_b_conditions_bis}) are satisfied, then
$$
W^T W \left(\begin{array}{c} a \\ - b \end{array}\right)
%= \left(\begin{array}{c c} 0 & \beta \\ \beta^T & 0 \end{array}\right) \left(\begin{array}{c} a \\ - b \end{array}\right)
= \left(\begin{array}{c} a - \beta b \\ \beta^T a - b \end{array}\right) = \left(\begin{array}{c} (2 - \lambda) a \\ (\lambda - 2) b \end{array}\right) = (2 - \lambda) \left(\begin{array}{c} a \\ - b \end{array}\right)
$$
and thus $(a,-b)$ is an eigenvector of $W^T W$ of eigenvalue $2 - \lambda$ that proves that all $2 n$ eigenvalues of $W^T W$ are in couples $\lambda_i, 2 - \lambda_i$ and since all eigenvalues are strictly positive we can conclude that $0 < \lambda_i < 2$ and thus we can write $\Lambda_{W^2} = \left(\begin{array}{c c} \Lambda & \\ & 2 - \Lambda \end{array}\right)$. Getting back to the solution of (\ref{formula_a_b_conditions_bis}) there are two different cases: if $\lambda = 1$ then (\ref{formula_a_b_conditions_bis}) reduces to
\begin{equation}
\label{formula_a_b_conditions_ter}
\left\{ \begin{array}{l l l}
\beta b & = & 0 \\
\beta^T a & = & 0 \dotinformula
\end{array} \right.
\end{equation}
that correspond to eigenvalues $0$ of $\beta$ and $\beta^T$ that are perfectly legal since, as we have seen the allowed eigenvalues of $t$ are in $\{ -1, e^{\pm i \theta} \}$ and consequently those of $\beta = \frac{1}{2} (\Identity + t)$ are in $\{ 0, \frac{1}{2} (1 + e^{\pm i \theta}) \}$. When $\lambda \ne 1$ we get easily that any $a \in \R^n$, part of an eigenvector of $W^T W$, must satisfy
\opt{margin_notes}{\mynote{mbh.note: alternatively $\beta^T \beta b = (\lambda - 1)^2 b$}}%
\begin{equation}
\label{formula_a_b_conditions_2_bis}
\beta \beta^T a = (\lambda - 1)^2 a
\end{equation}
and thus we can conclude that in both cases $a$ is an eigenvector of $\beta \beta^T$ of eigenvalue $(\lambda - 1)^2 \ge 0$. Since the $n \times n$ matrix $\beta \beta^T$ is symmetric it has $n$ orthogonal eigenvectors $a_i \in \R^n$ that we put in an orthogonal matrix $A$ and thus
$$
A^T \beta \beta^T A = \Lambda_\beta^2 \qquad \implies \qquad \beta \beta^T A = A \Lambda_\beta^2
$$
where $\Lambda_\beta^2$ is an $n \times n$ diagonal matrix with $(\lambda - 1)^2$ on the diagonal. With (\ref{formula_a_b_conditions_2_bis}) we can write the eigenvector of $W^T W$ associated to eigenvalues $\lambda_i \ne 1$
$$
\left(\begin{array}{c c} \Identity & \beta \\ \beta^T & \Identity \end{array}\right) \left(\begin{array}{c} a_i \\ (\lambda_i - 1)^{-1} \beta^T a_i \end{array}\right) = \left(\begin{array}{c} a_i + (\lambda_i - 1)^{-1} \beta \beta^T a_i \\ \beta^T a_i + (\lambda_i - 1)^{-1} \beta^T a_i \end{array}\right) = \left(\begin{array}{c} \lambda_i a_i \\ \lambda_i (\lambda_i - 1)^{-1} \beta^T a_i \end{array}\right)
$$
%and supposing that there are $k \le n$ of these eigenvectors we can arrange them in a $2 n \times k$ matrix of the form
%$$
%\left(\begin{array}{c} a_i \\ (\lambda_i - 2)^{-1} \beta^T a_i \end{array}\right) \qquad i = 1,2, \ldots, k \dotinformula
%$$
while for the remaining eigenvectors of $W^T W$ associated to eigenvalues $\lambda_j = 1$ of $W^T W$ and eigenvalues $0$ of $\beta \beta^T$ we get from (\ref{formula_a_b_conditions_ter})
$$
\left(\begin{array}{c c} \Identity & \beta \\ \beta^T & \Identity \end{array}\right) \left(\begin{array}{c} a_j \\ b_j \end{array}\right) = \left(\begin{array}{c} a_j \\ b_j \end{array}\right)
$$
where $a_j$ belongs to the null space of $\beta^T$ and $b_j$ to that of $\beta$ (\ref{formula_a_b_conditions_ter}). We define a $n \times n$ matrix $B$ formed by vectors $(\lambda_i - 1)^{-1} \beta^T a_i$ for eigenvectors $a_i$ corresponding to $\lambda_i \ne 1$ and by vectors $b_j$ for the cases of $\lambda_j = 1$ and it is an easy task to verify that also $B$ is orthogonal since $(\lambda_i - 1)^{-1} a_i^T \beta b_j = 0$.
%
%with which we can complete the set of $n$ eigenvectors of $W^T W$ in an $2 n \times n$ matrix we indicate with $\left(\begin{array}{c} A \\ B \end{array}\right)$ where $A$ is the orthogonal matrix of the eigenvectors of $\beta \beta^T$ and $B$ is built as seen .
%
%; from this relation we get immediately that the matrix $\beta^T A \Lambda_\beta^{-1}$, where clearly $\Lambda_\beta^{-1} \Lambda_\beta = \Identity$, is an orthogonal matrix and thus an element of $\OO{n}$.
%
In summary we can write the matrix $P$ of the $2 n$ eigenvectors of $W^T W$
$$
P
%= \left(\begin{array}{c c} A & A \\ \beta^T A \Lambda^{-1} & - \beta^T A \Lambda^{-1} \end{array}\right) :
= \frac{1}{\sqrt{2}} \left(\begin{array}{c c} A & A \\ B & - B \end{array}\right)
$$
with $A$ and $B$ both orthogonal and in $\OO{n}$ and thus
$$
P^T P = \frac{1}{2} \left(\begin{array}{c c} A^T & B^T \\ A^T & - B^T \end{array}\right) \left(\begin{array}{c c} A & A \\ B & - B \end{array}\right) = \left(\begin{array}{c c} \Identity & 0 \\ 0 & \Identity \end{array}\right)
$$
%$$
%R^T R = \left(\begin{array}{c c} A^T & \Lambda^{-1} A^T \beta \\ A^T & - \Lambda^{-1} A^T \end{array}\right) \left(\begin{array}{c c} A & A \\ \beta^T A \Lambda^{-1} & - \beta^T A \Lambda^{-1} \end{array}\right) =
%$$
that confirms that $P$ is the orthogonal rotation of $\OO{2 n}$ that diagonalizes $W^T W$; by Corollary~\ref{Witt_is_on_basis} $P$ is also an orthogonal Witt basis of $\R^{n,n}$.
\end{proof}

Two other props

\begin{MS_Proposition}
\label{MTNPs_O(2)}
In $\R^{2,2}$ given any 2 proper rotations (or antirotations) $t_\theta, t_\phi \in \OO{2}$ then for any $\phi \ne \theta$ for their corresponding MTNP
\opt{margin_notes}{\mynote{mbh.note: notation $t_\theta$ should be properly defined...}}%
$$
(x, t_\theta (x)) \cap (x, t_\phi (x)) = \{ 0 \} \quad \quad (x, t_\theta (x)) \oplus (x, t_\phi (x)) = \R^{2,2} \dotinformula
$$
\end{MS_Proposition}
\begin{proof}
Without loss of generality we assume $\theta = 0$ and $t_0 = \Identity$, a proper rotation, thus $(x, t_0 (x)) = P$ (\ref{formula_P_Q_def}) so we just need to prove that, for any $\phi \ne 0$ then $P \cap (x, t_\phi (x)) = \{ 0 \}$ and this follows immediately observing that by Proposition~\ref{condition_Witt_basis} $\det \left(\begin{array}{r r} \Identity & \Identity \\ \Identity & t_\phi \end{array}\right) = \sin^2 \frac{\phi}{2}$ that is not zero for $0 < \phi < 2 \pi$. Moving to antirotations without loss of generality we assume that for $\theta = 0$ $t_0 = \left(\begin{array}{r r} 1 & 0 \\ 0 & -1 \end{array}\right)$ and thus $(x, t_0 (x)) = \my_span{p_1, q_2}$ so we just need to prove that, for any $\phi \ne 0$, $\my_span{p_1, q_2} \cap (x, t_\phi (x)) = \{ 0 \}$ and the proof is identical to previous case.
\opt{margin_notes}{\mynote{mbh.note: and $\det () = - \sin^2 \frac{\phi}{2}$}}%
\end{proof}

\begin{MS_Proposition}
\label{MTNPs}
In $\R^{n,n}$ there are exactly $2^n$ different MTNP's that correspond to the $2^n$ elements of the Fock basis $\myFockB$ of its spinor space.
\end{MS_Proposition}

\begin{proof}
By definition of \mysetM{} we know that in $\R^{n,n}$ there are the $2^n$ MTNP of the Fock basis $\myFockB$ of spinor space.
\opt{margin_notes}{\mynote{mbh.note: commented here there is a longer proof by induction}}%
%
%We exploit bijection (\ref{formula_bijection_Nn_O(n)}) between MTNP's of \mysnqG{} and $\OO{n}$ and proceed by induction on $n$, for $n=1$ the Proposition is true since $\mbox{O}(1) = \{ \pm 1\}$ and the unique Witt basis of $\R^{1,1}$ is given by $\left(\begin{array}{r r} 1 & 1 \\ 1 & -1 \end{array}\right)$ called in physics the light cone in the $(x,t)$ plane. This basis uses both elements of $\mbox{O}(1)$ and thus all MTNP's of $\R^{1,1}$. Let the Proposition be true for $n-1$ and let us add two new dimension of generators $\mygen_{n}, \mygen_{2 n}$ orthogonal to previous linear space $\R^{n-1,n-1}$ and with which we can form two new null vectors $1/2 (\mygen_{n} \pm \mygen_{2 n})$. By induction hypothesis $\R^{n-1,n-1}$ contains $2^{n-1}$ different MTNP's and to any of these subspaces we can add one or the other of the two new null vectors so that we can build $2^{n}$ different MTNP's.
It remains to be proved that there are no other MTNP's that do not appear in \mysetM{} but this descends from Proposition~\ref{prop_2_BudinichP_1989} that shows that given $P$ any other MTNP's of \mysnqG{}, in an appropriate basis of $\R^{n,n}$, is one of the $2^n$ MTNP's of \mysetM{} and of the Fock basis $\myFockB$ of spinor space.
\end{proof}

%\newpage
\subsubsection{Isometries as elements of the algebra (first part of an old version with inner elements)}
\label{isometries_inner_old}
Any Clifford algebra contains its Clifford group that, in turn, contains the orthogonal group of the linear space, $\OO{n,n}$ in our case \cite[Chapter~13]{Porteous_1981}. Any $t \in \OO{n,n}$ is represented by an algebra element $T \in \myClg{}{}{\R^{n,n}}$ whose action on $x \in \R^{n,n}$ is given by $T x T^{-1}$ and the group product is given by the Clifford product of the corresponding algebra elements.
\opt{margin_notes}{\mynote{mbh.note: commented there is the introduction of $T$ as inner isomorphisms (less general and not needed in our case)}}%
%
%Since the dimension of the linear space $\R^{n,n}$ is even all automorphisms of $\myClg{}{}{\R^{n,n}}$ are inner and thus to any isometry $t \in \OO{n}$ associated to MTNP $(x, t(x))$ corresponds $T \in \myClg{}{}{\R^{n,n}}$ and the corresponding action on any element $x \in \myClg{}{}{\R^{n,n}}$ is given by $T x T^{-1}$ and this holds trivially also if $x$ is an element of the linear space $\R^{n,n}$.
%%
%\opt{margin_notes}{\mynote{mbh.note: probably I could simplify this stuff introducing the Clifford group and its subgroups $\OO{n,n}$ etc.}}%
%%

It is simple to verify that for example to the rotation that inverts $\mygen_{2 i}$ and $\mygen_{2 j}$ corresponds $\Lambda = \mygen_{2 i} \mygen_{2 j}$ and to the antirotation that inverts $\mygen_{2 i}$ corresponds $\Lambda = \omega \mygen_{2 i}$, the Hodge dual of $\mygen_{2 i}$. In summary to the isometries $\lambda \in \O1n$ correspond inner elements $\Lambda$ of the form
\begin{equation}
\label{formula_Lambda_def}
\Lambda = \omega^s \mygen_{2 i_1} \mygen_{2 i_2} \cdots \mygen_{2 i_s}
\end{equation}
\opt{margin_notes}{\mynote{mbh.ref: see pp. 58, 59}}%
that invert $s$ timelike generators, the factor $\omega^s$ accounting for antirotations when $s$ is odd. So for example in $\myClg{}{}{\R^{2,2}}$ it is simple to check that $\Lambda = \omega \mygen_{3} \in \O12$ and
$$
\Lambda P \Lambda^{-1} = \omega \mygen_3 \, \my_span{p_1, p_2} \, \mygen_3^{-1} \omega^{-1} = \my_span{q_1, p_2} \dotinformula
$$

We are now ready to define the isometry $\Lambda_j \in \O1n$ associated to a clause $z_j = \literal_{j_1} \literal_{j_2} \cdots \literal_{j_k}$ with $k$ literals; let $P_j = \my_span{p_{j_1}, p_{j_2}, \ldots, p_{j_k}}$ be the corresponding $k$TNP subspace of $P$ and let $r$ of the $k$ literals of $z_j$ appear in plain form and $k-r$ in complementary form, then $M(z_j)$ is defined as in (\ref{M_z_j_def}). Let
\begin{equation}
\label{formula_Lambda_j_def_old}
\Lambda_j = \omega^r \mygen_{2 j_1} \cdots \mygen_{2 j_r}
\end{equation}
with the understanding that $\Lambda_j = \Identity$ if $r = 0$, then it is simple to verify that $M(z_j) = \Lambda_j P_j \Lambda_j^{-1}$ and thus
\begin{equation}
\label{formula_Lambda_j_action_old}
M(z_j) \subset \Lambda_j P \Lambda_j^{-1}
\end{equation}
namely that the isometry $\Lambda_j$ (\ref{formula_Lambda_j_def_old}) transforms $P$ in a MTNP of \mysetM{} that has $M(z_j)$ as a subspace: in practice $\Lambda_j$ exchanges $r$ of the $n$ $p_i$ of $P$ with $q_i$. The definition of $\Lambda_j$ (\ref{formula_Lambda_j_def_old}) satisfying (\ref{formula_Lambda_j_action_old}) is not unique since we can freely add to $\Lambda_j$ (\ref{formula_Lambda_j_def_old}) any subset of the $n-k$ generators $\mygen_{2 i}$ whose indexes do not appear in $z_j$ and (\ref{formula_Lambda_j_action_old}) continues to hold. There are thus $2^{n-k}$ possible definitions of $\Lambda_j$ all satisfying (\ref{formula_Lambda_j_action_old}) and their action on $P$ produce the $2^{n-k}$ different MTNP $M(\psi_{z_j})$ of Corollary~\ref{M_z_j_coro} corresponding to the EFB expansion (\ref{formula_psi_z_j_EFB}).

We remark that for any choice $\Lambda_j \in \O1n$ and that all $\Lambda_j$ are involutions since $\Lambda_j^2 = \pm \Identity$.

We define the set of isometries induced by a clause $z_j$ as
\begin{equation}
\label{formula_cal_T_j_def_old}
{\cal T}_j' := \{\Lambda_j \in \O1n : M(z_j) \subset \Lambda_j P \Lambda_j^{-1} \}
\end{equation}
namely the set of $2^{n-k}$ elements of \O1n that satisfy (\ref{formula_Lambda_j_action_old}) and we can prove.... (omissis)

\subsubsection*{Things to do}
\begin{itemize}
\item ...
\end{itemize}

%\newpage
\newpage

\subsection*{$\R^{n,n}$ and its bases}
The whole story is narrated in, non degenerate, ``neutral'' linear space $\R^{n,n}$ (the familiar Minkowski plane being $\R^{1,1}$); $\R^{n,n}$ has a standard o.n. basis $E = \{\mygen_{1}, \mygen_{2}, \ldots, \mygen_{n}, \mygen_{n+1}, \ldots, \mygen_{2 n} \}$ (b.t.w. $\mygen_{i}$ are also the generators of Clifford algebra $\myClg{}{}{\R^{n,n}}$) and the generic scalar product of two vectors of this space $v_1, v_2$ may be written as (I use here standard $\cdot$ notation for scalar product to avoid introducing anticommutators that we don't really need afterwards)
$$
v_1 \cdot v_2 = v_1^T \left(\begin{array}{r r} \Identity & 0 \\ 0 & - \Identity \end{array}\right) v_2
$$
whereas the same scalar product in Euclidean space $\R^{2 n}$ would be just $v_1^T v_2$. The o.n. property of basis $E$ is
$$
E^T \left(\begin{array}{r r} \Identity & 0 \\ 0 & - \Identity \end{array}\right) E = \left(\begin{array}{r r} \Identity & 0 \\ 0 & - \Identity \end{array}\right)
$$
trivial in we assume $E = \Identity_{2 n}$.
For $\R^{n,n}$ we can always define the Witt, or null, basis:
\begin{equation}
\label{formula_LFR_Witt_basis}
\left\{ \begin{array}{l l l}
p_{i} & = & \frac{1}{2} \left( \mygen_{i} + \mygen_{i + n} \right) \\
q_{i} & = & \frac{1}{2} \left( \mygen_{i} - \mygen_{i + n} \right)
\end{array} \right.
\end{equation}
whose orthogonality properties are
$$
p_i \cdot p_j = q_i \cdot q_j = 0 \qquad p_i \cdot q_j = \delta_{i j} \qquad \forall i,j = 1, \ldots, n
$$
that are simple to check. Forming with the $2 n$ null vectors $p_i$ and $q_i$ two $2n \times n$ matrix $P$ and $Q$, we can condensate these relations in block matrix format as
\begin{equation}
\label{formula_LFR_2}
(P, Q)^T \left(\begin{array}{r r} \Identity & 0 \\ 0 & - \Identity \end{array}\right) (P, Q) = \left(\begin{array}{r r} 0 & \Identity \\ \Identity & 0 \end{array}\right) \dotinformula
\end{equation}

This stuff is very simple and can be found in many books of linear algebra and is usually called Witt decomposition of a linear space; frequently is more complicated because in spaces $\R^{m,n}$ with $m \ne n$ things are not so neat and, on top of that - God bless the mathematicians - one can consider degenerate spaces. My favourite book for this stuff is \cite[Chapter~9]{Porteous_1981}.

\subsection*{Null subspaces of $\R^{n,n}$}
Clearly $\my_span{p_1, p_2, \ldots, p_n}$ is a totally null subspace of $\R^{n,n}$ of maximum dimension $n$, or Maximally Totally Null Plane, MTNP for short. Stretching the notation we write $P = \my_span{p_1, p_2, \ldots, p_n}$ and similarly for $Q$, it is simple to verify that $\R^{n,n} = P \oplus Q$, the Witt decomposition of $\R^{n,n}$.

\opt{margin_notes}{\mynote{what follows is essentially p. 3 of my paper}}%
Let \mysnqG{} be the set of all MTNP of $\R^{n,n}$, \mysnqG{} is isomorphic to the subgroup $\OO{n}$ of $\OO{n,n}$ ({\bf this is ``Crucial Property'' {\#}1}). We can easily understand this relation: seeing the neutral space $\R^{n,n}$ as $\R^n \times \R^n$ we can write its generic element as $(x,y)$ and then (here by $x^2$ we mean $x \cdot x$, how is true in Clifford algebra)
$$
(x,y)^2 = x^2 - y^2
$$
and so for any $x \in \R^n \times \{0\}$ and $t \in \OO{n}$ $(x,t(x))$ is a null vector since
$$
(x,t(x))^2 = x^2 - t(x)^2 = x^2 - x^2 = 0
$$
and thus as $x$ spans the (spacelike) subspace $\R^n \times \{0\}$ then $(x,t(x))$ spans a MTNP of $\R^{n,n}$; we indicate this MTNP with self-explanatory notation as $(\Identity, t)$. Isometry $t \in \OO{n}$ establishes the quoted isomorphism since any MTNP of \mysnqG{} can be written in the form $(\Identity, t)$ \cite[Corollary~12.15]{Porteous_1981} and thus the one to one correspondence between MTNP of $\R^{n,n}$ and $t \in \OO{n}$ is formally established by the bijection
\begin{equation}
\label{formula_LFR_bijection_Nn_O(n)}
\OO{n} \to \mysnqG{}; t \to (\Identity, t) \qquad \implies \qquad \mysnqG = \{ (\Identity, t) : t \in \OO{n} \} \dotinformula
\end{equation}
\noindent For example in this setting the form of two generic null vectors of $P$ and $Q$ are respectively $(x, x)$ and $(x, -x)$ and in our new notation we represent the whole MTNP $P$ and $Q$ with
\begin{equation}
\label{formula_LFR_P_Q_def}
\begin{array}{l}
P = (\Identity, \Identity) \\
Q = (\Identity, -\Identity)
\end{array}
\end{equation}
and with this notation we may rewrite (\ref{formula_LFR_2}), duly normalized, as
\begin{equation}
\label{formula_LFR_3}
\frac{1}{2} (P, Q)^T \left(\begin{array}{r r} \Identity & 0 \\ 0 & - \Identity \end{array}\right) (P, Q) = \frac{1}{2} \left(\begin{array}{r r} \Identity & \Identity \\ \Identity & - \Identity \end{array}\right)^T \left(\begin{array}{r r} \Identity & 0 \\ 0 & - \Identity \end{array}\right) \left(\begin{array}{r r} \Identity & \Identity \\ \Identity & - \Identity \end{array}\right) = \left(\begin{array}{r r} 0 & \Identity \\ \Identity & 0 \end{array}\right) \dotinformula
\end{equation}

Also this stuff is simple and up to now I found it only in my favourite book \cite[Chapter~12]{Porteous_1981} but, as before, the book contains much much more than what is resumed here.

\subsection*{Simple spinors and null subspaces of $\R^{n,n}$}
Cartan was the first to prove that there exists a one to one correspondence between a subset of spinors, the \emph{simple} spinors, and a subset of all MTNP \mysnqG{}. Before looking into this relation we remark that simple spinors, that are $2^n$ for $\myClg{}{}{\R^{n,n}}$, are also linearly independent and thus form a natural basis for spinor space $S$: the Fock basis $\myFockB$.

How this basis can be formed is succintly explained in the paper at page 2 and, more in detail, in several of the quoted references. In a nutshell, given the null vectors $p_i$ and $q_i$ of the Witt basis of $\R^{n,n}$, we can form $2^n$ ``objects'' (multivectors actually) choosing one null vector, $p_i$ or $q_i$, for any $i = 1, \ldots, n$; for example in $\myClg{}{}{\R^{3,3}}$ $q_1 q_2 q_3$ is one of the $8$ simple spinors that form the Fock basis $\myFockB$ of its spinor space $S$.

The relevant property for us ({\bf ``Crucial Property'' {\#}2}) is the one to one correspondence between one of the $2^n$ simple spinors $\psi$ and the $2^n$ MTNP $M(\psi)$ of \mysnqG{} that can be formed choosing one element for each of the $n$ couples $(p_i, q_i)$. In previous example of $\myClg{}{}{\R^{3,3}}$ to simple spinor $\psi = q_1 q_2 q_3$ correspond the MTNP $M(\psi) = \my_span{q_1, q_2, q_3}$ and this is proved observing that for any $v \in M(\psi)$ then
$$
v \psi = (\alpha_1 q_1 + \alpha_2 q_2 + \alpha_3 q_3) q_1 q_2 q_3 = 0
$$
how is easily verified since $q_i^2 = 0$ and orthogonal vectors anticommute in Clifford algebra.

This one to one correspondence defines the set of MTNP $\mysetM = \{ M(\psi) \}$, namely the set of $2^n$ MTNP associated to the set of $2^n$ simple spinors of the Fock basis $\myFockB$; \mysetM{} is a subset of the larger set \mysnqG{} (\ref{formula_LFR_bijection_Nn_O(n)}).

This stuff is not-so-trivial and rarely tackled in Clifford algebras: up to now I found it only in my father's paper \cite{BudinichP_1989} and in Cartan book.

\subsection*{Simple spinors and SAT problems}
Given $n$ Boolean variables $x_i$ they can form $2^n$ different combinations (from to (T,T, \ldots, T) to (F,F, \ldots, F)) and any logical formula - and a fortiori any SAT problem - made with these $n$ Boolean variables is either T or F on any of the possible $2^n$ different combinations of the Boolean variables, more technically Boolean \emph{atoms}. A SAT problem is unSATisfiable if and only if its logical formula is F for \emph{all} the $2^n$ Boolean atoms. All this is more or less trivial.

A quite general property states that ``in any associative, unital, algebra - like $\myClg{}{}{\R^{n,n}}$ - every family of commuting, orthogonal, idempotents generates a Boolean algebra'' ({\bf ``Crucial Property'' {\#}3}). Here, with some technicalities, we see that a family of commuting, orthogonal, idempotents is given by $2^n$ simple spinors of $\myClg{}{}{\R^{n,n}}$ and thus, in turn, by the $2^n$ MTNP of $\mysetM \subset \mysnqG$.

We will thus embed any Boolean formula $S$ (not to be mistaken with spinor space $S$) in $\myClg{}{}{\R^{n,n}}$ and evaluate $S$ on Boolean atoms, indicated with $\rho_1, \rho_2, \ldots, \rho_n$ where $\rho_i \in \{ x_i, \myconjugate{x}_i \}$: we write this as $S(\rho_1, \rho_2, \ldots, \rho_n) \in \{T, F\}$.

In $\myClg{}{}{\R^{n,n}}$ we can calculate this value as $S \psi$ where $\psi = \rho_1 \rho_2 \cdots \rho_n$ is the simple spinor associated to the given Boolean atom: in summary by this embedding of $S$ in Clifford algebra ({\bf ``Crucial Property'' {\#}4})
$$
S(\rho_1, \rho_2, \ldots, \rho_n) \equiv F \qquad \iff \qquad S \psi = 0
$$
and thus SAT problems take an algebraic form in Clifford algebra.

In this section we have been running a little bit fast; this stuff appears only in my papers: this one and the previous one \cite{Budinich_2017}.

\subsection*{Wrapup {\#}1}
We put together various results we met up to now: in any Boolean problem formulated in Clifford algebra the following are fully equivalent:
\begin{itemize}
\item one of the $2^n$ possible assignments of the $n$ Boolean variables $\rho_1, \rho_2, \ldots, \rho_n$,
\item the simple spinor $\psi = \rho_1 \rho_2 \cdots \rho_n$,
\item the MTNP $M(\psi) \in \mysetM \subset \mysnqG$,
\item $t \in \OO{n}$ corresponding by (\ref{formula_LFR_bijection_Nn_O(n)}) to $M(\psi)$,
\item an orthogonal $n \times n$ matrix of $\R(n)$ corresponding to $t \in \OO{n}$
\end{itemize}
and from now on we will switch freely between these forms.

\subsection*{Formulation of SAT problems in Clifford algebra}
\opt{margin_notes}{\mynote{this part is just a resumè of pp. 5--8 of my paper}}%
We wrote an expression for an assignment of $n$ Boolean variables $\rho_1, \rho_2, \ldots, \rho_n$ in Clifford algebra but we still lack the formulation of a given SAT problem $S$ in this setting. This part appears in Section~3 of my paper or, more in detail, in my previous paper \cite{Budinich_2017}. In a nutshell any SAT problem is the logical AND of $m$ clauses; any clause being made by the logical OR of $k$ literals
$$
{\cal C}_j \equiv (\literal_{j_1} \lor \literal_{j_2} \lor \cdots \lor \literal_{j_k})
$$
that by De Morgan's relations can be written (remember that in Clifford algebra the logical AND is replaced by Clifford product)
$$
{\cal C}_j \equiv \myconjugate{\myconjugate{\literal}_{j_1} \myconjugate{\literal}_{j_2} \cdots \myconjugate{\literal}_{j_k}} : \equiv \Identity - z_j \qquad \mbox{with} \qquad z_j = \myconjugate{\literal}_{j_1} \myconjugate{\literal}_{j_2} \cdots \myconjugate{\literal}_{j_k}
$$
and $z_j$ is just the assignment of $k < n$ literals. With opportune conditions (Proposition~\ref{M_z_j} of my paper) this assignment correspond to a Totally Null Plane TNP of (non maximal) dimension $k < n$ indicated with $M(z_j)$.

At this point I report Propositions~\ref{SAT_in_TNP}, \ref{isomorphism_restricted} and \ref{SAT_in_n_O(1)} of my paper slightly adapted to present formalism
\begin{MS_Proposition}
\label{LFR_SAT_in_TNP}
A \SAT{} problem $S$ is unsatisfiable if and only if, for any of the $2^n$ assignments $\literal_1 \literal_2 \cdots \literal_n$ there exists at least one clause $z_j$ such that $M(z_j) \subseteq M(\literal_1 \literal_2 \cdots \literal_n)$ ({\bf ``Crucial Property'' {\#}5}).
\end{MS_Proposition}
Given considerations of Wrapup {\#}1 it is not surprising that this can be transposed to elements of $\OO{n}$: for this step we introduce the abelian group $\mbox{O}(1) \times \mbox{O}(1) \cdots \times \mbox{O}(1) = \stackrel{n}{\times} \mbox{O}(1) := \O1n$ that is a subgroup of $\OO{n}$. In familiar representation of elements of $\OO{n}$ with the $n \times n$ matrices of $\R(n)$, $\O1n$ is the abelian subgroup of diagonal matrices $\lambda \in \R(n)$ with $\pm 1$ on the diagonal, then
\begin{MS_Proposition}
\label{LFR_isomorphism_restricted}
The isomorphism (\ref{formula_LFR_bijection_Nn_O(n)}) when restricted to the subgroup \O1n of $\OO{n}$ has for image $\mysetM \subset \mysnqG$ ({\bf ``Crucial Property'' {\#}6}).
\end{MS_Proposition}
\noindent and Proposition~\ref{LFR_SAT_in_TNP} rephrased in terms of orthogonal transformations is:
\begin{MS_Proposition}
\label{LFR_SAT_in_n_O(1)}
A given \SAT{} problem $S$ with m clauses $z_j$ is unsatisfiable if and only if the isometries induced by its clauses (\ref{formula_cal_T_j_def}) form a cover for \O1n:
\begin{equation}
\label{formula_LFR_SAT_in_n_O(1)}
\cup_{j = 1}^m {\cal T}_j' = \O1n \dotinformula
\end{equation}
\end{MS_Proposition}

\subsection*{MTNP in basis independent form}
We have seen that any MTNP can be represented by $(\Identity, t)$ with $t \in \OO{n}$ but in general $t$ can take different forms depending on $\R^{n,n}$ basis.

The final step is to rephrase Proposition~\ref{LFR_SAT_in_TNP}, that is a relation on incidences between TNP subspaces of $\R^{n,n}$, in basis independent form since the property of intersection of subspaces cannot depend on the choice of a particular basis.

This is done by Theorem~\ref{SAT_in_O(n)_old} of my paper that generalizes Proposition~\ref{LFR_SAT_in_n_O(1)} to a basis independent form and shows that in this case the isometries associated with clauses form a cover, not only of \O1n, but of the full parent group $\OO{n}$ ({\bf ``Crucial Property'' {\#}7}).

The proof is based on a truly remarkable property of MTNP of $\R^{n,n}$ that essentially says: in $\R^{n,n}$ there are \emph{only} the $2^n$ MTNP of $\mysetM$ that ultimately correspond to $\lambda \in \O1n$, all other $t \in \OO{n}$ forming \mysnqG{} are just one of the $2^n$ MTNP in another basis; in other words in basis independent form $\mysetM = \mysnqG$.

Also this is a rarely cited property that appears in my father's paper \cite{BudinichP_1989} as Proposition~2 and, again, in Cartan book.

\subsection*{Wrapup {\#}2}
The difference between Proposition~\ref{LFR_SAT_in_n_O(1)} and Theorem~\ref{SAT_in_O(n)_old} of my paper is that \O1n{} is a \emph{discrete} group and to check if a set of isometries form a cover one has to check essentially all $2^n$ elements of \O1n{} to see if they are in the cover. Conversely $\OO{n}$ is a smooth, compact, disconnected real manifold of dimension $n (n-1)/2$ (the size of orthogonal matrices of $\R(n)$).

Would $\OO{n}$ be a flat space of dimension $n (n-1)/2$ (and it is not) we could simply take the set of isometries associated to the $m$ clauses of a SAT problem to see if they form a basis of this space, for example calculating their determinant, and in the affirmative we could deduce that they cover the full $\OO{n}$ and thus necessarily its subgroup \O1n{} and that the problem at hand is unSATisfiable.

This argument is flawed because $\OO{n}$ is not flat: in general it correspond to the sphere of dimension $n-1$, for example for $\OO{2}$ is the circle described by one parameter: the angle of rotation $\theta$ etc. etc.

} % note finali: stampate solo se all'inizio c'è l'opzione final_notes

\end{document}